\documentclass[reqno]{amsart}
\usepackage{amssymb}
\usepackage{amscd}
\usepackage{amsmath}

\usepackage{color}  
\definecolor{darkgreen}{rgb}{0.03, 0.5, 0.03}
 \newcommand{\ec}{\color{black}} %

 \newcommand{\X} {\boldsymbol X}
 \newcommand{\co} {\boldsymbol c}
 \newcommand{\Na} {\mbox{\rm Na}}
  \newcommand{\Max} {\mbox{\rm Max}}
    \newcommand{\Spec} {\mbox{\rm Spec}}
      \newcommand{\QMax} {\mbox{\rm QMax}}
        \newcommand{\QSpec} {\mbox{\rm QSpec}}
\newcommand{\hgt} {\mbox{\rm ht}}

\newcommand{\fbar}{\boldsymbol{\overline{F}}}

  \newcommand{\stf} {\star{_{\!{_f}}}}
    \newcommand{\stt} {\widetilde{\star}}
    
      \newcommand{\cl}{^{{\mbox{\rm \footnotesize{cl}}^{{\star}}}}}
	  
	    \newcommand{\clt}{^{{\mbox{\rm 
	  \footnotesize{cl}}^{{\widetilde{\star}}}}}}

\newtheorem{theorem}{Theorem}[section]
\newtheorem{lemma}[theorem]{Lemma}

\newtheorem{example}[theorem]{Example}
\newtheorem{remark}[theorem]{Remark}

\newtheorem{corollary}[theorem]{Corollary}

\begin{document}

\title[Uppers to zero and Pr\"ufer-like domains]
{Uppers to zero in polynomial rings and Pr\"ufer-like domains}

\author[G. W. Chang and M. Fontana]{Gyu Whan Chang  \; \;\;\; \;\;\;\;\;\; \;\;\; \ \;\;\; \;\;\; \;\;\;\;\;\; \;\;\; \;\;\;  Marco Fontana \\
Department of Mathematics \;\;\; \;\;\; \;\;\; Dipartimento di Matematica \\
University of Incheon, Korea \;\;\;\;\;  \; \;\; Universit\`a ``Roma Tre'', Italy}

\date{}

%

 \thanks{\it Email.  \rm G.W.C.:  \texttt{whan@incheon.ac.kr}; \; M.F.: \texttt{fontana@mat.uniroma3.it}}%



\begin{abstract}

 Let $D$ be an integral domain and $X$ an indeterminate over $D$. 
It is well known that  (a) $D$ is quasi-Pr\"ufer (i.e, its integral closure is a Pr\"ufer domain)  if and only if  each upper to zero $Q$  in   $D[X] $  contains a polynomial $g \in D[X]$
 with  content   $\co_D(g) = D$; (b) an upper to zero $Q$ in $D[X]$ is a maximal
$t$-ideal if and only if $Q$ contains a    nonzero   polynomial $g \in D[X]$
 with  $\co_D(g)^v = D$. 
Using these facts,  the notions of UM$t$-domain (i.e., an integral domain such that  each upper to zero is a maximal
$t$-ideal) and quasi-Pr\"ufer domain can be naturally  extended to  the semistar operation setting and studied in a unified frame.
In this paper, given a semistar operation $\star$ in the sense of Okabe-Matsuda, we introduce 
the $\star$-quasi-Pr\"ufer domains. We give several characterizations of these domains and we investigate  their relations with the UM$t$-domains and the Pr\"ufer $v$-multiplication domains.   
\end{abstract}

\maketitle


\section*{introduction  and background results}

 Gilmer and Hoffmann characterized Pr\"ufer
domains as those integrally closed domains $D$, such that the extension of $D$ inside its quotient field  is a primitive extension \cite[Theorem 2]{GH}. (Relevant definitions and results are reviewed in the sequel.)  Primitive extensions are strictly related with relevant properties of the prime spectrum of the polynomial ring. In particular, from the previous characterization it  follows that a Pr\"ufer domain is an integrally closed   quasi-Pr\"ufer  domain (i.e., an integral domain such that each prime ideal of the polynomial ring  contained in an extended  prime is extended \cite{ace}) \cite[Section 6.5]{fhp}.  A quasi-Pr\"ufer domain $D$ can be characterized by the fact that  each upper to zero $Q$   in   $D[X] $  contains a polynomial $g \in D[X]$
with content $\co_D(g) = D$ (Theorem  \ref{qp}).   On the other hand, a ``weaker'' version of the last property can be used for characterizing   upper to zero that are maximal $t$-ideals in the polynomial ring.    Recall that $D$ is called a 
UM$t$-domain  (UM$t$ means ``upper to zero is a maximal $t$-ideal") if
every upper to zero in $D[X]$ is a maximal $t$-ideal   \cite[Section 3]{hz} and this happens if and only if    each upper to zero in $D[X]$ contains a nonzero polynomial  $g \in
D[X]$ with   $\co_D(g)^v = D$ \cite[Theorem 1.1]{fgh}. 
Using the previous observations,  the notions of UM$t$-domain and quasi-Pr\"ufer domain can be naturally  extended  to   the semistar operation setting and studied in a unified frame.
More precisely, given a semistar operation $\star$ in the sense of Okabe-Matsuda \cite{o-m}, we introduce in a natural way 
the $\star$-quasi-Pr\"ufer domains    and semistar analog of other relevant notions like  primitive extension and  incomparability (INC) property .  We give  several characterizations of   the $\star$-quasi-Pr\"ufer domains  and we investigate   their relations with the UM$t$-domains and the Pr\"ufer $v$-multiplication domains  \cite{griffin}. 

 More precisely, let  $\star$ be a semistar operation on an integral domain $D$ with quotient field $K$. Among other things, we prove that $D$ is a  $\stf$-quasi-Pr\"ufer domain if and only if   
$ D \subseteq K$ is a $\stf$-primitive extension, if and only if $D$ is a  $\stf$-INC-domain,  
if and only if  each  overring $R$ of $D$ is a   $(\stf)_{\iota}$-quasi-Pr\"ufer domain, 
where $\iota: D \hookrightarrow R$ is the canonical embedding, if and only if
every prime ideal of $\Na(D, \star_f)$ is extended from $D$, if and only if
$\Na(D, \stf)$ is a quasi-Pr\"ufer domain, if and only if the integral closure of
$\Na(D, \stf)$   is a Pr\"ufer domain, if and only if $D_P$ is a quasi-Pr\"ufer domain, 
for each quasi-$\stf$-maximal ideal $P$  of $D$. Moreover, we show that if $\star$ is a (semi)star operation, then $D$ is a $\stf$-quasi-Pr\"ufer domain if and only if $D$ is  a $t$-quasi-Pr\"ufer domain  and each  $\stf$-maximal ideal  of $D$ is a $t$-ideal (equivalently, $\widetilde{\ \stf} = w$).

We  also   show that this general approach sheds new light on some delicate aspects of the classical theories.
In particular,  we give a contribution to the open problem of whether the integral closure
of a UM$t$-domain is a P$v$MD by showing that $D$ is a UM$t$-domain if and only if  the  $w$-closure $\widetilde{D}$ of $D$ is a P$v$MD and  the $w$-operations on $D$ and  $\widetilde{D}$ are related by $(w_D)_{\widetilde{\iota}} =w_{\widetilde{D}}$, where $\widetilde{\iota}:D \hookrightarrow \widetilde{D}$ is the canonical embedding.  Moreover,  among other results,  we provide a   positive answer to a Zafrullah's conjecture on the local-global behaviour of the UM$t$ domains  \cite[page 452]{zafrullah}.  
$$
\ast \:  \ast \:  \ast \:  \ast \:  \ast \: 
$$

\ec Let $D$ be an integral domain with quotient field $K$. Let
$\boldsymbol{\overline{F}}(D)$ denote the set of all nonzero
$D$--submodules of $K$ and let $\boldsymbol{F}(D)$ be the set of
all nonzero fractional ideals of $D$, i.e. $E \in
\boldsymbol{F}(D)$ if $E \in \boldsymbol{\overline{F}}(D)$ and
there exists a nonzero $d \in D$ with $dE \subseteq D$. Let
$\boldsymbol{f}(D)$ be the set of all nonzero finitely generated
$D$--submodules of $K$. Then, obviously $\boldsymbol{f}(D)
\subseteq \boldsymbol{F}(D) \subseteq
\boldsymbol{\overline{F}}(D)$.

Following Okabe-Matsuda \cite{o-m}, a \emph{semistar operation} on $D$ is a map $\star:
\boldsymbol{\overline{F}}(D) \to \boldsymbol{\overline{F}}(D), E
\mapsto E^\star$,  such that, for all $x \in K$, $x \neq 0$, and
for all $E,F \in \boldsymbol{\overline{F}}(D)$, the following
properties hold:
\begin{enumerate}
\item[$(\star_1)$] $(xE)^\star=xE^\star$;
 \item[$(\star_2)$] $E
\subseteq F$ implies $E^\star \subseteq F^\star$;
\item[$(\star_3)$] $E \subseteq E^\star$ and $E^{\star \star} :=
\left(E^\star \right)^\star=E^\star$.
\end{enumerate}

 Recall that,  given a
semistar operation $\star$ on $D$, for all $E,F \in
\boldsymbol{\overline{F}}(D)$,  the following basic formulas
follow easily from the axioms:    $$\begin{array}{rl} (EF)^\star
=& \hskip -7pt (E^\star F)^\star =\left(EF^\star\right)^\star
=\left(E^\star
F^\star\right)^\star\,;\\
(E+F)^\star =& \hskip -7pt \left(E^\star + F\right)^\star= \left(E
+
F^\star\right)^\star= \left(E^\star + F^\star\right)^\star\,;\\
(E:F)^\star \subseteq & \hskip -7pt (E^\star :F^\star) = (E^\star
:F) =
\left(E^\star :F\right)^\star,\;\, \mbox{\rm if \ } (E:F) \neq 0;\\
(E\cap F)^\star \subseteq & \hskip -7pt E^\star \cap F^\star =
\left(E^\star \cap F^\star \right)^\star,\;\, \mbox{\rm if \ }
E\cap F \neq (0)\,;
\end{array}
$$
 \noindent cf.  for instance \cite[Theorem 1.2 and p.  174]{FH2000}.

A \emph{(semi)star operation} is a semistar operation that,
restricted to $\boldsymbol{F}(D)$,  is a star operation (in the
sense of \cite[Section 32]{gilmer}). It is easy to see that a
semistar operation $\star$ on $D$ is a (semi)star operation if and
only if $D^\star = D$.

 If $\star$ is a semistar operation on $D$, then we can
consider a map\ $\star_{\!_f}: \boldsymbol{\overline{F}}(D) \to
\boldsymbol{\overline{F}}(D)$ defined, for each $E \in
\boldsymbol{\overline{F}}(D)$, as follows:

\centerline{$E^{\star_{\!_f}}:=\bigcup \{F^\star\mid \ F \in
\boldsymbol{f}(D) \mbox{ and } F \subseteq E\}$.}

\noindent It is easy to see that $\star_{\!_f}$ is a semistar
operation on $D$, called \emph{the semistar operation of finite
type associated to $\star$}.  Note that, for each $F \in
\boldsymbol{f}(D)$, $F^\star=F^{\star_{\!_f}}$.  A semistar
operation $\star$ is called a \emph{semistar operation of finite
type} if $\star=\star_{\!_f}$.  It is easy to see that
$(\star_{\!_f}\!)_{\!_f}=\star_{\!_f}$ (that is, $\star_{\!_f}$ is
of finite type).

If $\star_1$ and $\star_2$ are two semistar operations on $D$, we
say that $\star_1 \leq \star_2$ if $E^{\star_1} \subseteq
E^{\star_2}$, for each $E \in \fbar(D)$. This is equivalent to say
that $\left(E^{\star_{1}}\right)^{\star_{2}} = E^{\star_2}=
\left(E^{\star_{2}}\right)^{\star_{1}}$, for each $E \in
\fbar(D)$.  Obviously, for each semistar operation $\star$ defined on $D$, we
have $\star_{\!_f} \leq \star$. Let $d_D$ (or, simply, $d$)  be the \it identity (semi)star operation on $D$,   \rm clearly $d \leq \star$, for all semistar operation $\star$ on $D$.

We say that a nonzero ideal $I$ of $D$ is a
\emph{quasi-$\star$-ideal} if $I^\star \cap D = I$, a
\emph{quasi-$\star$-prime} if it is a prime quasi-$\star$-ideal,
and a \emph{quasi-$\star$-maximal} if it is maximal in the set of
all   proper   quasi-$\star$-ideals. A quasi-$\star$-maximal ideal is  a
prime ideal. It is possible  to prove that each  proper   quasi-$\star_{_{\!
f}}$-ideal is contained in a quasi-$\star_{_{\! f}}$-maximal
ideal.  More details can be found in \cite[page 4781]{fl}. We
will denote by $\QMax^{\star}(D)$  (resp., $\QSpec^\star(D)$) the set of the 
quasi-$\star$-maximal ideals  (resp., quasi-$\star$-prime ideals) of $D$.
When $\star$ is a (semi)star operation the notion of  quasi-$\star$-ideal coincides with the ``classical'' notion of  \it  $\star$-ideal \rm (i.e., a nonzero ideal $I$ such that $I^\star = I$). 

The
$\star$-dimension of $D$, denoted by $\dim^\star(D)$, is defined by the
supremum of $\{n \mid  P_1 \subsetneq P_2 \subsetneq \cdots \subsetneq P_n$ is a
chain of  quasi-$\star$-prime ideals of $D\}$. Thus,  when $\star$ is a semistar operation of finite type and $D$ is not a field,    $\dim^\star(D)=1$ if and
only if each quasi-$\star$-maximal ideal of $D$ has height-one. 

  If $\Delta$ is a set of prime ideals of an integral domain
  $D$,  then the semistar operation $\star_\Delta$ defined
  on $D$ as follows

  \centerline{$
  E^{\star_\Delta} := \bigcap \{ED_P \;|\;\, P \in \Delta\}\,,
  \;  \textrm {  for each}    \; E \in \boldsymbol{\overline{F}}(D)\,,
  $}

  \noindent is called \it the spectral semistar operation associated to
  \rm
  $\Delta$.
  A semistar operation $\star$ of an integral domain $D$ is
  called
  \it a
  spectral semistar operation \rm if there exists a subset $
  \Delta$ of the prime spectrum of $D$, $\mbox{\rm Spec}(D)$,  such that $\,\star =
  \star_\Delta\,$.

  When $\Delta := \QMax^{\star_{_{\! f}}}(D)$, we set $\stt:= \star_{\Delta}$, i.e.

  \centerline{  
  $E^{\stt} := \bigcap \left \{ED_P \mid P \in  \QMax^{\star_{_{\! f}}}(D) \right\}$,  \;  for each $E \in \boldsymbol{\overline{F}}(D)$.
  }

 A semistar operation $\star$ is \emph{stable} if $(E \cap F)^\star
= E^\star \cap F^\star$, for each $E,F \in \fbar(D)$.
Spectral semistar operations are stable \cite[Lemma~4.1 (3)]{FH2000}. In particular, $\stt$ is a semistar operation stable and of finite type \cite[Corollary 3.9]{FH2000}.

 By $v_D$ (or, simply, by $v$) we denote  the $v$--(semi)star
operation defined as usual by  $E^v := (D:(D:E))$, for each $E\in
\boldsymbol{\overline{F}}(D)$. By  $t_D$ (or, simply, by $t$) we
denote  $(v_D)_{_{\! f}}$ the $t$--(semi)star operation on $D$ and
by  $w_D$ (or just by $w$) the stable semistar operation of finite
type associated to $v_D$ (or, equivalently, to $t_D$), considered
by F.G.  Wang  and R.L. McCasland in \cite{WMc97} (cf. also \cite{gv});  i.e.,  $w_D :=
\widetilde{v_D} = \widetilde{t_D}$.  Clearly $w_D\leq t_D \leq v_D$.  Moreover, it is easy to see that for each   (semi)star operation $\star$ of $D$, we have $\star \leq v_D$ and $\stf \leq t_D$ (cf. also \cite[Theorem 34.1 (4)]{gilmer}).

 If $I \in \boldsymbol{\overline{F}}(D)$, we say that $I$ is
\emph{$\star$--finite} if there exists $J \in \boldsymbol{f}(D)$
such that $J^\star=I^\star$.  It is immediate to see that if
${\star_1} \leq {\star_2}$ are semistar operations and $I$ is
${\star_1}$--finite, then $I$ is ${\star_2}$--finite. In
particular, if $I$ is $\star_{\!_f}$--finite, then it is
$\star$--finite. The converse  is not true and it is possible to
prove that $I$ is $\star_{\!_f}$--finite if and only if there
exists $J \in \boldsymbol{f}(D)$, $J \subseteq I$, such that
$J^\star=I^\star$  \cite[Lemma 2.3]{fo-pi}.

If $I$ is a nonzero ideal of $D$, we say that $I$ is
\emph{$\star$--invertible} if $(II^{-1})^\star = D^\star$.   From the
definitions and from the fact that $\QMax^{\stf}(D) = \QMax^{\stt}(D)$ \cite[Corollary   3.5(2)]{fl} it follows easily   that an ideal $I$ is  
$\tilde{\star}$--invertible if and only if $I$ is $\star_{_{\!
f}}$--invertible.  If $I$ is $\star_{_{\! f}}$--invertible, then
$I$ and $I^{-1}$ are $\star_{_{\! f}}$--finite \cite[Proposition
2.6]{fo-pi}.

 ~Let $R$ be an overring of an integral domain $D$, let $\iota: D \hookrightarrow R$ be the canonical embedding 
and  let $\star$ be a semistar operation of $D$. We denote by  $\star_\iota$   the semistar operation of $R$ defined by $E^{\star_\iota} := E^\star$,  for each $E \in \fbar(R) \ (\subseteq \fbar(D))$.   Let $\ast$ be a semistar operation of $R$ and let $\ast^\iota$ be the semistar operation on $D$ defined by $E^{\ast^\iota} := (ER)^\ast$,  for each $E \in \fbar(D)$.    It is not difficult to see that 
${(\ast^{\iota})}{_{_{\! f}}} = ({\ast}{_{_{\!f}}})^{\iota}$  and  if $\star$ is a semistar operation of finite type (resp., a stable semistar operation) of $D$ then ${\star_\iota}$ is a semistar operation of finite type (resp., a stable semistar operation) of $R$ (cf. for instance \cite[Proposition 2.8]{fl2} and \cite[Propositions 2.11 and 2.13]{pi}).


\section{quasi-pr\"ufer domains}

Let $D$ be an integral domain  with quotient field $K$,  and let  ${\boldsymbol X}$  be a nonempty set
of indeterminates over $K$.  For each polynomial $f \in K[\X]$, we denote by $\co_D(f)$ (or, simply, $\co(f)$) \it the content on $D$ of the polynomial $f$, \rm i.e., the  (fractional)  ideal of $D$ generated by the coefficients of $f$. For  each fractional ideal $J$ of $D[\X]$,  with $J \subseteq K[\X]$,   we denote by $\co_D(J)$  (or, simply, $\co(J)$) the (fractional) ideal $\{\co_D(f) \mid f \in J\}$ of $D$. Obviously, for each ideal $J$ in  $D[\X]$,  $J\cap D \subseteq \co_D(J)$ and $(J\cap D)[\X] \subseteq J \subseteq \co_D(J)[\X] $.

Taking the properties of prime ideals in polynomial extensions of Pr\"ufer domains as a starting point, the quasi-Pr\"ufer notion was introduced in \cite{ace} for arbitrary rings (not necessarily domains).     As in \cite[page  212]{fhp}, we
say that $D$ is a {\em quasi-Pr\"ufer domain} if for each prime
ideal $P$ of $D$, if $Q$ is a prime ideal of $D[{\bf \X}]$ with $Q
\subseteq P[{\bf \X}]$, then  $Q = (Q \cap D)[{\bf \X}]$.  It is well known that an integral domain is a Pr\"ufer domain if and only if it is integrally closed and quasi-Pr\"ufer \cite[Theorem 19.15]{gilmer}.

Consider now the following condition:
\begin{center} ({\bf qP$^\prime$}) if $Q$ is a prime ideal of
$D[{\bf \X}]$ with  $\co_D(Q) \subsetneq D$,  then $Q = (Q \cap D)[{\bf
\X}]$. \end{center} 

It is clear that $D$ satisfies ({\bf qP$^\prime$}) if and only if $D$ is a
quasi-Pr\"ufer domain.  Therefore,  an integrally closed
domain $D$ is a Pr\"ufer domain  if and only if $D$ satisfies ({\bf qP$^\prime$}).

 Let $D \subseteq R$ be an extension of integral domains, and let
$P$ be a prime ideal of $D$. We say that \it $D \subseteq R$ satisfies
INC at $P$ \rm  if whenever $Q_1$ and $Q_2$ are prime ideals of $R$
such that $Q_1 \cap D = P = Q_2 \cap D$, then $Q_1$ and $Q_2$ are
incomparable.  If  $D \subseteq R$ satisfies
INC at every prime ideal of $D$, $D \subseteq R$ is said \it an INC-extension. \rm  The domain $D$ is \it an INC-domain \rm  if, for each overring $R$ of $D$, $D\subseteq R$ is an INC-extension.

An element $u \in R$ will be said to be {\em primitive over $D$} if $u$ is a root
of  a primitive polynomial on $D$ (i.e., a nonzero polynomial $f \in D[X]$ with $\co_D(f) = D$). The extension $D\subseteq R$ is called \it a primitive extension \rm  (or, \it a P-extension \rm \cite{GH})  if each element of $R$ is primitive over $D$.  

  A nonzero prime ideal $Q$ in the polynomial ring $D[X]$ is called \it an upper to zero \rm (McAdam's terminology) if $Q \cap D =(0)$. Let $P:= Q \cap D$,  if $Q = P[X]$ then   $Q$ is called \it an extended prime \rm  of $D[X]$ (more details can be found in \cite{hous}).  

Recall that Gilmer and Hoffmann characterized Pr\"ufer
domains as those integrally closed domains $D$, such that the   embedding    of $D$ inside its quotient field  is a P-extension \cite[Theorem 2]{GH},   and that D. Dobbs  in \cite{dobbs} characterized  P-extensions in terms of INC-domains.    The natural link between quasi-Pr\"ufer domains and  primitive extensions 
is recalled in the following theorem, where we collect several useful characterizations of quasi-Pr\"ufer domains   (cf. also the  very recent survey paper by E. Houston \cite{hous}).


\begin{theorem} \label{qp}
Let $D$ be an integral domain with quotient field $K$,   let $X$ be an indeterminate over $D$ and let   $\mathcal N := \{ g \in D[X] \mid \co_D(g) = D \}$  be the set of primitive polynomials over $D$.   Then the following statements are
equivalent.
\begin{enumerate}
\item $D$ is a quasi-Pr\"ufer domain.
\item[(1$'$)] $D$ satisfies ({\bf qP$^\prime$}) for one indeterminate.
\item [(2)] Each upper to zero  in $D[X]$ contains a polynomial  $g \in
D[X]$ with $\co_D(g)  = D$.
\item[(3)] If $Q$ is an upper to zero in $D[X]$, then $ \co_D(Q)  = D$.
\item[(4)] $D \subseteq K$ is a primitive extension.
\item[(5)] $D$ is an INC-domain.
\item[(6)] The integral closure of $D$ is a Pr\"ufer domain.
\item[(7)] Each overring of $D$ is a quasi-Pr\"ufer domain.
\item[(8)] Each prime ideal of  $D[X]_{\mathcal N}$   is extended from $D$.

\item[(9)] $D[X]_{\mathcal N}$  is a quasi-Pr\"ufer domain.
\item [(10)]The integral closure of  $D[X]_{\mathcal N}$  is a Pr\"ufer
domain.
\item[(11)] $D_M$ is a quasi-Pr\"ufer domain,  for each maximal ideal $M$ of $D$.
 
\end{enumerate}
\end{theorem}

\begin{proof}
  (1)$\Leftrightarrow$(4)$\Leftrightarrow$(5)$\Leftrightarrow$(6)$\Leftrightarrow$(7)  and 
 (9)$\Leftrightarrow$(10) by
\cite[Corollary 6.5.14]{fhp}. Moreover,  (3)$\Leftrightarrow$(6)  by \cite[Theorem 2.7]{adf}  

  (2)$\Leftrightarrow$(3), (1)$\Leftrightarrow$(11) and  (1)$\Rightarrow\! (1')$   are clear.   

  $(1')\! \Rightarrow$(3)   If  $Q$ is an upper to zero, then $Q \neq (Q\cap D)[X]$, 
 and   thus, by $(1')$, $\co_D(Q) = D$.  

 (6)$\Leftrightarrow$(10)   Let $\overline{D}$ be the integral closure
of $D$, and let $\overline{{\mathcal N}}  :=   \{h \in \overline{D}[X] \mid  \co_{\overline{D}}(h)  =
\overline{D}\}$. Then  it is clear that   $\overline{D}[X]_{{\mathcal N}} = \overline{D}[X]_{\overline{{\mathcal N}}}$.   Moreover,   
$\overline{D}[X]_{\mathcal N}$  coincides with   the integral closure of $D[X]_{\mathcal N}$  \cite[Chapter 5, Proposition 5.12 and Exercise 9]{am}. Finally, recall that  $\overline{D}$
is a Pr\"ufer domain if and only if  $\overline{D}[X]_{\overline{{\mathcal N}}}$   is a Pr\"ufer
domain \cite[Theorem 33.4]{gilmer}.

 $(1')\! \Rightarrow$(8)   Let $ \boldsymbol{\mathfrak{Q}}$ be a prime ideal of $D[X]_{\mathcal N}$. Then
$\boldsymbol{\mathfrak{Q}} = QD[X]_{\mathcal N}$ for some prime ideal $Q$ of $D[X]$. Since $ \boldsymbol{\mathfrak{Q}}
\subsetneq D[X]_{\mathcal N}$, $Q\cap {\mathcal N} = \emptyset$; hence $\co(Q) \subsetneq
D$. So, by $(1')$,
$Q = (Q \cap D)[X]$.  
Thus $ \boldsymbol{\mathfrak{Q}}= {(Q \cap D)}D[X]_{\mathcal N}$.

 (8)$\Rightarrow\! (1')$ Let $Q$ be prime ideal of $D[X]$ such that $\co(Q) \subsetneq
D$. Then $Q \cap {\mathcal N} = \emptyset$ and thus $QD[X]_{\mathcal N}$ is a prime ideal of $D[X]_{\mathcal N}$. Therefore, by (8),  $(Q\cap D)D[X]_{\mathcal N}=QD[X]_{\mathcal N}$, and 
 hence  $Q = (Q\cap D)[X]$.  
\end{proof}

In view of the extensions to the case of semistar operations, we introduce the following notation. Let $\star$ be a semistar operation on $D$,  if     $\mathcal N^\star
:=  \{ g\in D[X] \mid  g \neq 0 \mbox{ and }  \co_D(g)^\star  =D^\star\}$,   then we set $\Na(D, \star) := D[X]_{\mathcal N^{\star}}$. The ring of rational functions $\Na(D, \star)$ is called \it  the $\star$--Nagata domain of $D$.  \rm When $\star =d$ the identity (semi)star operation on $D$,  ${\mathcal N}^d ={\mathcal N}$ (the multiplicative set of $D[X]$ introduced  in Theorem \ref{qp}) and we set simply $\Na(D)$ instead of $\Na(D, d)= D[X]_{\mathcal N}$.  Note that  $\Na(D)$ coincides with the classical Nagata domain $D(X)$ (cf. for instance \cite[Chapter I, \S 6 page 18]{n} and \cite[Section 33]{gilmer}).

Recall  from \cite[Propositions 3.1 and 3.4]{fl} that:
\begin{enumerate}
\item[(a)]  $\mathcal N^\star= \mathcal N^{\stf}= \mathcal N^{\stt}=  D[X] \setminus \bigcup \{P[X] \mid P
  \in \QMax^{\stf}(D) \}$
is a saturated multiplicatively closed subset of  $D[X]$.

 \item[(b)]    $\Na(D, \star) = \Na(D, \stf) =\Na(D, \stt) =\bigcap\{ D_P(X) \mid P
  \in \QMax^{\stf}(D) \}$.

 \item[(c)]   
  $ \QMax^{\stf}(D)= \{M \cap D \mid M \in  \Max(\Na(D, \star)) \}\,.$

  \item[(d)]   $E^{\tilde{\star}} = E\!\cdot\!\Na(D,\star)
  \cap K$,   for each $E \in
  {\overline{\boldsymbol{F}}}(D)$.
  
  \end{enumerate}


\begin{remark}\label{umt}  \rm  
 (a)  It is well known that an upper to zero $Q$ in $D[X]$ is a maximal
$t$-ideal if and only if $Q$ contains a nonzero polynomial   $g \in D[X]$
with $\co_D(g)^t \ (= \co_D(g)^v ) =D$    \cite[Theorem 1.4]{hz}.     Recall that $D$ is called a {\em
UM$\ \!\!t$-domain}   if
every upper to zero in $D[X]$ is a maximal $t$-ideal   \cite[Section 3]{hz}.    An overring $R$ of $D$ is called \it $t$-linked to $D$ \rm  if, for each nonzero finitely generated ideal $I$ of $D$, $(D:I) =D$ implies $(R:IR) = R$  (cf. for instance  \cite{dhlz} and \cite{dhlz2}).   Recall that UM$t$-domains can be characterized by weaker ($t$--) versions of some of the statements of Theorem \ref{qp}, since  \it the following statements are equivalent:
\begin{enumerate} 
\item[$(1_{t})$ ]  $D$ is a UM\ \!$t$-domain. 
 
\item[$(2_{t})$ ]  Each upper to zero in $D[X]$ contains a   nonzero  polynomial  $g \in
D[X]$ with   $\co_D(g)^t = \co_D(g)^v  = D$. 

\item[$(3_{t})$ ]   If $Q$ is an upper to zero in $D[X]$, then  $ \co_D(Q)^t  = D$.  

\item[$(7_{t})$ ]   Each $t$-linked overring to $D$ is a UM\ \!$t$-domain. 
 
\item[$(8_{t})$ ] Each prime ideal of $\Na(D, t)$ is extended from $D$.

\item[$(11_{t})$ ]  $D_P$ is a quasi-Pr\"ufer domain,  for each maximal $t$-ideal $P$ of $D$.
\end{enumerate} \rm  
 For the proof see \cite[Theorem 1.1]{fgh}  and \cite[Theorem 2.6 (1)$\Leftrightarrow$(8)]{cz}.  

(b)  Note that if $P \subseteq Q$ are two primes ideals in a UM$t$-domain with $P \neq (0)$ and if $Q$ is a  prime $t$-ideal   then $P$ is also a   prime $t$-ideal   \cite[Corollary 1.6]{fgh}.

(c) With the notation introduced just before this remark, one of the arguments in the proof of (6)$\Leftrightarrow$(10) in Theorem \ref{qp} shows that, for any integral domain $D$,  the integral closure of $\Na(D)$ is $\Na(\overline{D})$.
 
(d)  Recall that an integral domain $D$ is called a \it Pr\"ufer $v$-multiplication domain \rm (for short, \it P$v$MD \rm) if each nonzero finitely generated ideal of $D$ is $t$-invertible or, equivalently, if   $(FF^{-1})^t =D$, for each $F \in \boldsymbol{f}(D)$  \cite{griffin}.
 It is known that a domain $D$ is an integrally closed domain and a UM$t$-domain  if and and only if $D$ is a P$v$MD \cite[Proposition 3.2]{hz}. 
    But  M. Zafrullah \cite[page 452]{zafrullah} mentioned a problem that seems to be still open:  \sl is the integral closure
of a UM\ \!$t$-domain a P$v$MD~\!\!?  \    \rm    We will give some contributions to this problem in the following  Corollaries \ref{p*md-um*} and  \ref{umt-closure}.   

 A related question is the following: \sl if the integral closure $\overline{D}$ of an integral domain $D$ is a P$v$MD what can be said about the UM\ \!$t$-ness of $D$~\!\!?  \  \rm 
An answer to this question was recently given by Chang and Zafrullah \cite[Remark 2.7]{cz} where they provide an example of a non-UM$t$ domain with the integral closure which is a P$v$MD.    
 
  \end{remark} 

 Using the notion of UM$t$-domain (recalled in the previous remark) we have further characterizations of a quasi-Pr\"ufer domain (cf. Theorem \ref{qp}).  

\begin{corollary} \label{qp-umt}
The following statements are equivalent for an integral domain
$D$. \begin{enumerate}
\item $D$ is a quasi-Pr\"ufer domain.
\item[(12)]    Each overring of $D$ is a UM$\ \!t$-domain.  

\item[(13)]  $D$ is a  UM $\!t$-domain  and each maximal ideal of $D$ is a
$t$-ideal.
\item[(14)]  $D$ is a UM $\!t$-domain  and $d =w$. 

\end{enumerate}
 In particular, in a quasi-Pr\"ufer domain every nonzero prime ideal is a $t$-ideal.  
\end{corollary}

\begin{proof}  (1)$\Leftrightarrow$(12)  by \cite[Corollary 3.11]{fgh} and Theorem \ref{qp} ((1)$\Leftrightarrow$(6)). 

 (1)$\Rightarrow$(13)   If $Q$ is an upper to zero in $D[X]$, then
$Q$ contains a nonzero  polynomial    $g \in D[X]$ with $ \co_D(g)   = D$   (Theorem
\ref{qp} ((1)$\Rightarrow$(2))).  Clearly   $ \co_D(g)^t   = D$,  and thus $D$ is a UM$t$--domain 
 (Remark \ref{umt} (a)  or \cite[Theorem 1.4]{hz}). Let $M$ be a maximal ideal of $D$. If $ M^t
 =D$, there is a polynomial  $0\neq h \in M[X]$ such that $\co_D(h)^t  = D$. 
 It is easy to see that $hD[X]_{M[X]} \cap D = (0)$.   In this situation, there exists an upper to
zero $Q'$  in $D[X]$  such that  $hD[X] \subseteq Q' \subseteq M[X]$  
\cite[Lemma 1.1 (b)]{dl}. Hence $Q'$ (and thus $M[X]$) contains a
 nonzero  polynomial   $g'$ with $\co_D(g')  = D$ by Theorem
 \ref{qp} ((1)$\Rightarrow$(2)),  thus   $D[X] =\co_D(g')D[X] \subseteq M[X]$, a contradiction.
 Therefore $M^t \subsetneq D$, hence   $M$ is a $t$-ideal.

 (13)$\Rightarrow$(1)   Let $Q$ be an upper to zero in $D[X]$.  Since we are assuming that $D$ is a UM$t$-domain,  then
$Q$ is a maximal $t$-ideal of $D[X]$,  and hence $Q$ contains a
polynomial  $0\neq g \in D[X]$ with $\co_D(g)^t =\co_D(g)^v  = D$   \cite[Theorem 1.4]{hz}.
 Furthermore, by assumption,  if $M$ is a maximal ideal of $D$, then  $\co_D(g)  \nsubseteq M$ since
$M$ is a $t$-ideal. Hence $ \co_D(g)  = D$,    and thus $D$ is a
quasi-Pr\"ufer domain by Theorem  \ref{qp} ((2)$\Rightarrow$(1)).

 (13)$\Rightarrow$(14)   Note that from (13) it follows easily that $\Max(D) = \Max^t(D) $.  
Thus $d = \widetilde{d} = \widetilde{t} = w$.

(14)$\Rightarrow$(13)  Under the present assumption $\Max(D) = \Max^w(D)$ and it is known that 
$\Max^w(D) = \Max^t(D)$ (cf. for instance \cite[Corollary 3.5 (2)]{fl}).

The last statement is an easy consequence of the  fact that a quasi-Pr\"ufer domain is a UM$t$-domain and of Remark \ref{umt} (b).  
\end{proof}

\begin{remark} \rm (a)  From the previous Corollary \ref{qp-umt} ((1)$\Leftrightarrow$(13)), we easily deduce that
the condition $(11_{t})$ in Remak \ref{umt}  (a), that characterizes the UM$t$-domains,  is equivalent to the following:
\begin{enumerate}

\item[$(11^\prime_{t})$]  \it $D_P$ is a UM\ \!$t_{D_P}$-domain and $PD_P$ is a maximal $t_{D_P}$-ideal of $D_P$,   for each maximal $t_D$-ideal $P$ of $D$.  \rm  
\end{enumerate}
\noindent(Cf. also \cite[Theorem 1.5]{fgh} and \cite[Theorem 3.13]{hous}.)  This result provides a positive answer to the following  Zafrullah's conjecture  \cite[page 452]{zafrullah}: \sl an integral domain $D$ is a UM$t_D$-domain if and only if $D_M$ is a UM\ \!$t_{D_M}$-domain,  for each maximal  ideal $M$ of $D$, and  $D$ is well behaved \rm  (i.e.,  a domain such that prime $t$-ideals of the domain extend to prime $t$-ideals in the rings of fractions of $D$). 

As a matter of fact, the ``only if part'',  on which was based the conjecture,  was already proved in \cite[Propositions 1.2 and 1.4]{fgh}; the ``if part''  follows from the equivalence of   $(11^\prime_{t})$ with $(1_t)$ of Remark  \ref{umt}  (a).  

(b) Note that the condition (12) in the previous Corollary \ref{qp-umt} can be stated in the following equivalent form:
\begin{enumerate}
\item[$(12')$] \it $D$ is  a   UM \!$t$-domain  and each overring of $D$ is $t$-linked to $D$. \rm 
\end{enumerate}
(Cf. \cite[Theorem 2.4]{dhlrz}.)

(c)  In relation with (14) of Corollary \ref{qp-umt}, we recall that the domains for which $d =w$ were introduced and studied in \cite{mim} under the name of \it DW-domains \rm  (cf. also \cite{pt} for further information on these domains). \rm  A DW-domain  $D$ can be characterized by the property that each overring $R$ of $D$ is $t$-linked to $D$ (cf. \cite[Theorem 2.6]{dhlrz}, \cite {dhlz2} and \cite[Proposition 2.2]{mim}).

\end{remark}

\begin{corollary}
Let $D$ be a quasi-Pr\"ufer domain.
Then $\dim(D)=   \dim^w(D) =   \dim^t(D) =  \dim(\Na(D))$.  
\end{corollary}

\begin{proof}
This follows because,  in the present situation, $d = w$,  every nonzero prime ideal of $D$ is a
$t$-ideal  (Corollary \ref{qp-umt})  and each prime ideal of   $\Na(D)$  is extended from  $D$
by Theorem  \ref{qp} ((1)$\Rightarrow$(8)).  \end{proof}

\begin{remark} \rm Note that even in the Pr\"ufer domain case, it might happen that $\dim(D) \ (=\dim^t(D)) \gneq \dim^v(D)$.  For instance take a nondiscrete valuation domain. In this case, the maximal ideal is not a $v$-ideal.
\end{remark}

\section{$\star$-quasi-pr\"ufer domains and  uppers to zero}


 Let $\star$ be a semistar operation on an integral domain $D$.  We want to introduce a semistar analog  to the notion of quasi-Pr\"ufer domain and to   the related notion of  UM$t$-domain. 
 
 We say that an integral domain $D$ is a \it $\star$-quasi-Pr\"ufer  domain \rm if  the following property holds: 
\begin{enumerate}
\item[($\star${\bf{qP}})] if $Q$ is a prime ideal in $D[X]$ and $Q \subseteq P[X]$, for some $P\in \QSpec^{\star}(D)$, then  $Q = (Q\cap D)[X]$. 
\end{enumerate}

It is clear from the definition that the $d$-quasi-Pr\"ufer  domains are exactly  the  quasi-Pr\"ufer  domains.

\begin{lemma}\label{*qp}  \ Let $\star$ be a semistar operation on an integral domain $D$.  The following statements are equivalent:
\begin{enumerate}
\item[(i)]  $D$ is a $\star$-quasi-Pr\"ufer  domain.
\item[(ii)]  Let $Q$ be  an upper
to zero  in $D[X]$,  then
 $\co_D(Q)\not\subseteq P$, for each $P \in \QSpec^\star(D)$.
 \item[(iii)]   Let $Q$ be  an upper
to zero  in $D[X]$, then 
 $Q \not\subseteq P[X]$,   for each $P \in \QSpec^\star(D)$.
  \item[(iv)]  $D_P$ is a quasi-Pr\"ufer domain, for each $P \in \QSpec^\star(D)$.
 \end{enumerate}
\end{lemma}

\begin{proof}  (i)$\Rightarrow$(iii) follows immediately from the definition.

(iii)$\Rightarrow$(ii)  If $Q$ is an upper to zero then by assumption  
 $Q \not\subseteq P[X]$, for all $P\in \QSpec^{\star}(D)$. Then  $\co(Q)\not\subseteq P$, for each $P \in \QSpec^\star(D)$, since $Q \subseteq \co_D(Q)[X]$.
 
(ii)$\Rightarrow$(i)  Assume that $Q$ is a prime ideal in $D[X]$ such that $(Q\cap D)[X] \subsetneq Q \subseteq P[X]$, for some $P \in \QMax^\star(D)$. Then we can find an upper to zero $Q_1$ in $D[X]$ such that $Q_1 \subseteq Q$  \cite[Theorem A]{dl}. Thus $\co_D(Q_1) \subseteq \co_D(Q) \subseteq P$, for some $P \in \QSpec^\star(D)$, and this contradicts the present hypothesis.

(i)$\Rightarrow$(iv) Let $P \in \QSpec^\star(D)$.  In order to show that $D_P$ is a quasi-Pr\"ufer domain, we prove   the  condition $(1^\prime)$ of Theorem \ref{qp}.  If $Q$ is a prime ideal of $D_P[X]$ with $\co_{D_P}(Q) \subsetneq D_P$, then $\co_{D_P}(Q) \subseteq PD_P$, and hence $Q \subseteq PD_P[X]$. So $Q \cap D[X] \subseteq P[X]$, and by (i) we have $Q \cap D[X] = (Q \cap D)[X]$. Hence $Q = (Q \cap D_P)[X]$. Thus $D_P$ is a quasi-Pr\"ufer domain.

(iv)$\Rightarrow$(i) Let $Q$ be a prime ideal of $D[X]$ such that $Q \subseteq P[X]$ for some $P \in \QSpec^\star(D)$. Then   $QD_P[X] \subseteq PD_P[X]$, and hence $QD_P[X] = 
(QD_P[X] \cap D_P)[X]$ by (iv). Thus $Q = (QD_P[X]  \cap D_P)[X] \cap D[X] = (Q \cap D)[X]$.  
\end{proof}

Since a quasi-$\star$-ideal is also a quasi-$\stf$-ideal, it is clear that $\stf$-quasi-Pr\"ufer implies $\star$-quasi-Pr\"ufer.  Recall that every  quasi-$\stf$-ideal is contained in a quasi-$\stf$-maximal ideal and each quasi-$\stf$-maximal ideal is a prime ideal \cite[Lemma 2.3]{fl}.  Therefore, the set $\QSpec^{\stf}(D)$ is always nonempty.  On the other hand $\QSpec^{\star}(D)$ can be empty  and in this case the notion of  $\star$-quasi-Pr\"ufer domain can be very weak.  

Note also that,  when $\star$ is a semistar operation of finite type,  in the condition ($\star${\bf{qP}}) and
in the properties  (ii), (iii), and (iv)  of the previous Lemma \ref{*qp} we can replace $\QSpec^\star(D)$ with $\QMax^\star(D)$, obtaining equivalent statements.

\begin{example} \label{nqp} \sl Example of a $\star$-quasi-Pr\"ufer domain which is not a $\stf$-quasi-Pr\"ufer domain.

\rm  Let $W$ be a 1-dimensional nondiscrete valuation domain with maximal ideal $N$ and residue field $k:=W/N$.  Let $Z$ be an indeterminate over $W$.  Passing to the Nagata's ring $V:=W(Z)$,  it is wellknown that $V$ is also a 1-dimensional nondiscrete valuation domain, with maximal ideal $M:= N(Z)$ and residue field $k(Z)$ (cf.  \cite[Theorem 33.4]{gilmer} and \cite[Theorem 14.1 and Corollary 15.2]{hu}). Let $\pi: V \rightarrow V/N = k(Z)$ be the canonical projection and let $D =\pi^{-1}(k)$.  Clearly, $D$ is an integrally closed 1-dimensional pseudo-valuation domain with maximal ideal $M$ and with associated valuation overring $V= (M:M)$ \cite[Theorem 2.10]{hh}. Note that $V$ has no divisorial primes, since $M$ is not finitely generated \cite[Exercise 12, page 431]{gilmer} and that the $t$-operation on a valuation domain coincides with $d$ the identity (semi)star operation.
Let  $\iota: D \hookrightarrow V$ be the canonical embedding and let $\star := (v_V)^\iota$ be the semistar operation on $D$ defined by $E^\star :=(EV)^{v_V}$, for each $E \in \overline{\boldsymbol{F}}(D)$. Note that $\star$ is not of finite type and more precisely it is not difficult to see that:
$$
\stf = {((v_V)^{\iota})}_{_{\! f}}=
{((v_V){_{_{\! f}}})}^{\iota} = {(t_V)}^{\iota} = {(d_V)}^{\iota} \; \; \; \mbox{ \cite[Proposition 2.13]{pi}}. 
$$
 Therefore
$E^{\stf} =EV$, for each $E \in \overline{\boldsymbol{F}}(D)$.  In particular, $M$ is a (quasi-)$\stf $-maximal ideal of $D$.  
Note that $D$ is not a $\stf$-quasi-Pr\"ufer domain since,   if $X$ is an indeterminate over $D$, $\dim(D[X])=3$,   because there    exists  an upper to zero $Q$ in $D[X]$ contained in $M[X]$ \cite[Theorem 2.5 and Remark 2.6]{hh2}.   On the other hand  $D$ is trivially a $\star$-quasi-Pr\"ufer domain  since $D$ does not possess quasi-$\star$-prime ideals, because $M^\star = (MV)^{v_V} = M^{v_V} = V$.
\end{example}

Because of the previous observations and Example \ref{nqp},  we consider with a special attention the case of $\stf$-quasi-Pr\"ufer domains.

\begin{lemma}\label{*fqp}  \ Let $\star$ be a semistar operation on an integral domain $D$. 
The following   statements   are equivalent:

\begin{enumerate} 

\item[($1_{\stf}$)] $D$ is a $\stf$-quasi-Pr\"ufer  domain.
\item[($2_{\stf}$)]  Each
upper to zero in $D[X]$ contains a nonzero  polynomial $g \in D[X]$ with $\co(g)^\star = D^\star$.
\item[($3_{\stf}$)]  If $Q$  is an  upper
to zero  in $D[X]$, then 
 $\co(Q)^{\stf}= D^{\star}$. 

\end{enumerate} 
\end{lemma}

\begin{proof}  $(1_{\stf})\Leftrightarrow (3_{\stf})$ follows from Lemma \ref{*qp}   because the property  
 $Q \not\subseteq P[X]$, for all $P\in \QMax^{\stf}(D)$ is equivalent to 
$\co_D(Q)^{\stf} = D^\star$ (since each proper  quasi-$\stf$-ideal is contained in a quasi-$\stf$-maximal).

$(3_{\stf})\Rightarrow (2_{\stf})$  is obvious.

 $(2_{\stf})\Rightarrow (1_{\stf})$ Let $Q$ be a prime ideal in $D[X]$ such that $Q \subseteq P[X]$, for some $P\in \QSpec^{\stf}(D)$. Assume $(Q\cap D)[X] \subsetneq Q$.  Then we can find an upper to zero $Q_1$ in $D[X]$ such that $Q_1 \subseteq Q$  \cite[Theorem A]{dl}. By assumption, there exists a  nonzero polynomial $ g \in Q_1$ such that $\co_D(g)^{\star} = D^\star$, hence in particular $\co_D(Q_1)^{\stf} = D^\star$ and so 
$\co_D(Q)^{\stf} = D^\star$.  This implies that $Q \not\subseteq P[X]$, for all $P\in \QSpec^{\stf}(D)$ and this contradicts the assumption.
\end{proof}

\begin{corollary} \label{<}  Let $\star, \star_1$ and $ \star_2$ be semistar operations on an integral domain $D$.
\begin{enumerate}
\item[(a)]  \it    Assume that  $\star_1 \leq \star_2$.  If $D$ is  a $\star_1$-quasi-Pr\"ufer domain  then  $D$ is a $\star_2$-quasi-Pr\"ufer domain. 
\item[(b)] \it  $D$ is a $t$-quasi-Pr\"ufer domain if and only if $D$ is a UM\ \!$t$-domain.

\item[(c)] \it  $D$ is a $\stf$-quasi-Pr\"ufer domain if and only if $D$ is a $\stt$-quasi-Pr\"ufer domain.
\end{enumerate}
\end{corollary}
\begin{proof}
(a) and (b) follow easily from Lemma \ref{*fqp}  ($(1_{\stf})\Leftrightarrow (2_{\stf})$) and from  Remark \ref{umt} (a). For (c) note also that $\co_D(g)^\star = D^\star$ if and only $\co_D(g) \not\subseteq P$ for all $P \in \QMax^{\stf}(D)$ and that   $\QMax^{\stf}(D)= \QMax^{\stt}(D)$  \cite[Lemma 2.3 (1) and Corollary 3.5 (2)]{fl}.
\end{proof}

\begin{remark} \rm For $\star =v$, we have observed in Corollary \ref{<} (b) that  the $t$-quasi-Pr\"ufer domains coincide with the UM$t$-domains, i.e., the domains such that each upper to zero in $D[X]$ is a maximal  $t_{D[X]}$-ideal.  There is no immediate extension to the semistar setting of the previous characterization, since  in the general case we do not have the possibility to work at the same time with a semistar  operation (like the $t$-operation) defined both on $D$ and on $D[X]$. 

 At this point it is natural to formulate the following   question.  
 
 \noindent   {\bf Question}: Given a semistar operation  of finite type  $\star$ on $D$, is it possible to define in a canonical way a semistar operation  of finite type  $\star_{\!_{D[X]}}$ on $D[X]$, such that $D$ is a $\star$-quasi-Pr\"ufer domain if and only if each upper to zero in $D[X]$ is a quasi-$\star_{\!_{D[X]}}$-maximal ideal~\!?  \footnote{Added in Proofs: This problem was solved by the authors in case of a stable semistar operation of finite type. The corresponding paper ``Uppers to zero and semistar operations
in polynomial rings'' is now published in Journal of Algebra \bf 318 \rm(2007) 484--493.}
 
 However, we want to mention that  Okabe and Matsuda \cite[Definition 2.10]{om} introduced a star-operation analog of the notion of
UM$t$-domain: given a star operation $\ast$ on an integral domain $D$, they call $D$ a \it $\ast$-UMT ring \rm  if each upper to zero contains a  nonzero   polynomial $g \in D[X]$ with $\co_D(g)^\ast = D$. This notion coincides with the notion of $\ast_{_{\!f}}$-quasi-Pr\"ufer domain introduced above in the more general setting of the semistar operations (Lemma \ref{*fqp}).

 \end{remark}

The next goal is to extend to the case of general $\stf$-quasi-Pr\"ufer domains the characterizations given in Theorem \ref{qp}.   For this purpose we need to extend some definitions to the semistar setting.

Let $D \subseteq R$ be an extension of integral domains and let $\star$ be a semistar  operation on $D$. We will
say that $R$ is a {\em $\star$-INC-extension} of $D$ if whenever $Q_1$
and $Q_2$ are  nonzero  prime ideals of $R$  such that $Q_1 \cap D = Q_2 \cap
D$ and $(Q_1 \cap D)^{\star }\subsetneq D^\star$,   
then $Q_1$ and $Q_2$ are
incomparable.  We also say that $D$ is a {\em $\star$-INC-domain} if
each overring of $D$ is a $\star$-INC-extension of $D$.  Moreover, we say that  
 an element $u \in R$
is {\em $\star$-primitive over $D$} if $u$ is a root
of a nonzero polynomial $g \in D[X]$ with $\co_D(g)^\star = D^\star$. 

 Note that the notion of $d$-primitive (respectively, $d$-INC) extension coincides with the ``classical'' notion of 
primitive (respectively, INC) extension.  It is obvious that the notions of  $\star$-primitive and $\stf$-primitive coincide, while $\stf$-INC-extension  implies $\star$-INC-extension. 
The converse is not true as it will be shown in the following example.

\begin{example} \label{inc} \sl Example of a $\star$-INC extension which is not a $\stf$-INC extension. \rm

  Let $D, V, M$ and $\star$ be  as in Example \ref{nqp}. It is easy to see that  $D$ is not a $\stf$-INC domain.  For instance, if $R:= \pi^{-1}(k[Z])$,  then  $M$ is a prime ideal also in $R$ and all the maximal ideals of $R$ and the prime (non maximal) ideal $M$ of $R$ have the same trace in $D$, that is $M$.  Since $M$ is a (quasi-)$\stf$-maximal ideal of $D$,  $D \hookrightarrow R$ is not a $\stf$-INC extension. On the other hand  $D$ is vacuously a $\star$-INC domain (the only nonzero prime of $D$ is $M$ and $M^\star \cap D = D$). 
\end{example}

\begin{lemma}  Let $D$ be an integral domain with quotient field $K$ 
and let $P$ be a prime ideal of $D$. 
For $u \in K$, $D \subseteq D[u]$ satisfies INC at $P$ 
if and only if there exists $0 \neq g  \in D[X]$ such
that  $\co_D(g) \nsubseteq P$  and $g(u)= 0$. 
\end{lemma}

\begin{proof}
Let $I$ be   the kernel of the canonical surjective homomorphim $D[X] \rightarrow D[u], \ X \mapsto u$. 
 It is known that  $D \subseteq D[u]  \cong D[X]/I $   satisfies INC at $P$ if and only if $\co_D(I)
\nsubseteq P$ \cite[Proposition 2.0]{pa}. Suppose $\co_D(I) \nsubseteq
P$. Choose $\ec a \in \co_D(I) \setminus P$.  Since $a \in \co_D(I)$,  then there exist   a finite family of polynomials  $f_1,
\dots , f_k \in I$ such that $a \in \co_D(f_1)+ \co_D(f_2) + \cdots + \co_D(f_k)$.
Let $ g := f_1 + X^{n_1+1}f_2 + X^{n_1+n_2+2}f_3 + \dots +
X^{n_1+n_2 + \dots + n_{k-1} + k-1}f_k$, where $n_i$ is the degree
of $f_i$. Then $g \in I$,   $ a \in \co_D(g)$,  and $g(u)=0$. Since $a \notin P$, then $\co_D(g) \nsubseteq P$.  
Conversely, if $g (u) = 0$, then $g \in I$, and hence  
$\co_D(g)  \nsubseteq P$   implies $\co_D(I) \nsubseteq P$. 
\end{proof}

  Recall from Remark \ref{umt} (a) that an overring $R$ of $D$ is called  $t$-linked to $D$ if for each nonzero finitely generated ideal $I$ of $D$, $(D:I) =D$ implies $(R:IR) = R$.

\begin{remark}  \label{t-linked} \rm The notion of $t$-linked overring can be characterized in several ways.  In particular, the following  statements   are equivalent   \cite[Proposition 2.1]{dhlz2}:  

\begin{enumerate}
\item[(i)] \it $R$ is a $t$-linked overring to $D$.

\item[(ii)] \it    For  each nonzero finitely generated ideal $I$ of $D$, $I^{t_D}=D$ implies $(IR)^{t_R}= R$.

\item[(iii)] \it   For  each prime (or maximal)  $t_R$-ideal  $Q$ of $R$, $(Q \cap D)^{t_D} \subsetneq D$.

\end{enumerate}  
\end{remark}

In case that $\star$ is a semistar operation on $D$, we need the following (relativized) extension of the notion of $t$-linkedness. We say that an overring $R$ of $D$ is \it $t$-linked to $(D, \star)$ \rm if, for each nonzero finitely generated ideal $I$ of $D$, $I^{\star}=D^\star$ implies $(IR)^{t_R}= R$  \cite[Section 3]{eBF}. Therefore  the notion of ``$R$   is $t$-linked to $(D, t_D)$'' coincides with  the ``classical'' notion of  ``$R$   is  $t$-linked to $D$''.

We collect in the following lemma some characterizations of the $t$-linkedness in the semistar setting.

\begin{lemma} \label{t-*-linked}  Let $\star$ be a semistar operation on an integral domain $D$ with quotient field $K$ and let $R$ be an overring of $D$. 
The following   statements  are equivalent:
\begin{enumerate}
\item[(i)]  $R$ is a $t$-linked overring to $(D, \star)$.
\item[(i$_{\!_f}$)] $R$ is a $t$-linked overring to $(D, \stf)$.
\item[($\widetilde{\mbox{\rm i}}$)]  $R$ is a $t$-linked overring to $(D, \stt)$.

\item[(ii)] For each nonzero  ideal $I$ of $D$, $I^{\stf}=D^\star$ implies $(IR)^{t_R}= R$.

\item[(iii)] For each prime (or maximal) $t_R$-ideal $Q$ of $R$, $(Q \cap D)^{\stf} \subsetneq D^\star$.

\item[(iv)] For each  proper  $t_R$-ideal $J$ of $R$, $(J \cap D)^{\stf} \subsetneq D^\star$.
\item[(v)] $R = R^{\stt} \ (= R\!\cdot\! \Na(D, \star) \cap K)$.

\end{enumerate}
\end{lemma}

\begin{proof}  (i)$\Leftrightarrow$(i$_{\!_f}$)$\Leftrightarrow$($\widetilde{\mbox{\rm i}})$  because, for a nonzero finitely generated ideal $I$ of $D$, $I^{\star}= D^\star$ is equivalent to say that $I \not\subseteq P$, for all $P \in \QMax^{\stf}(D) = \QMax^{\stt}(D)$.

The equivalences (i)$\Leftrightarrow$(ii)$\Leftrightarrow$(iii)$\Leftrightarrow$(iv) are consequences of \cite[Proposition 3.2]{eBF}.

 (iii)$\Rightarrow$(v)  From the assumption it follows that,  for each  maximal $t$-ideal $Q$ of $R$, there exists a quasi-$\stf$-maximal ideal $P$ of $D$ containing $Q\cap D$ and thus $ D_P \subseteq  R_{D\setminus P} \subseteq R_Q$. Therefore  $R \subseteq R^{\stt} =  \bigcap \{RD_P \mid P \in \QMax^{\stf}(D) \} \subseteq  \bigcap \{R_Q \mid Q \in \Max^{t_R}(R) \}  = R$.  

 (v)$\Rightarrow$(iii) Let $Q$ be a prime $t$-deal of $R$ such that $(Q \cap D)^{\stf} = D^\star$. Therefore there exists a nonzero finitely generated ideal $I \subset Q\cap D$ such that $I^{\stf} = D^\star$. In particular, we have  $IR[X] \cap \mathcal{N}^{\stf} \neq \emptyset$ and so $(IR)^{\stt}= IR\!\cdot\!\Na(D, \star) \cap K = IR[X]_{{\mathcal{N}^{\stf} }} \cap K= R[X]_{{\mathcal{N}^{\stf} }} \cap K=   R^{\stt} =R$. On the other hand $IR \subseteq Q$ and so $ (IR)^{t_R} \subseteq Q^{t_R} =Q$. Moreover, if we denote by $\iota$  the canonical embedding of $D$ into $R$, then $\ast:=(\stt)_\iota$  is a (semi)star operation of finite type on $R$, since $R =R^{\stt}= R^\ast$.  Therefore $\ast \leq t_R$ and so we get a contradiction, since $ R =  (IR)^{\stt} = (IR)^\ast \subseteq (IR)^{t_R} \subseteq Q \subsetneq R$.
\end{proof}  

\begin{remark} \rm Given  a star operation $\ast$ on $D$,  the property (v) of Lemma  \ref{t-*-linked}  is used in  \cite[page 224]{chang2} for giving the definition ``$R$ is $\ast$-linked to $D$''  (teminology used in  that paper). That notion coincides with the notion  of ``$R$ is $t$-linked to $(D, \ast)$'' (terminology used here)  (cf. \cite[Proposition 3.2]{chang2}).  

Note also that, from the previous Lemma \ref{t-*-linked}, we re-obtain in particular the equivalences stated in Remark \ref{t-linked}.
\end{remark}

As a consequence of the previous   Lemma \ref{t-*-linked}  we deduce immediately the following two corollaries.

\begin{corollary}  Let $R$ be an overring of  an integral domain $D$ with quotient field $K$. Then $R$ is a $t$-linked overring to $D$ if and only if $R = R^{w_D} ( \ = R \cdot \Na(D, v_D) \cap K)$.
\end{corollary}

  For the next statement we need to recall the notion of $\star$-valuation overring
(a notion due essentially to {P. Jaffard} 
(\cite[page 46]{J}, see  also   \cite[Chapters 15 and 18]{HK2}).
 For a domain
$D$ and a semistar operation $\star$ on $D$,  we say 
that 
a
valuation overring $V$ of $D$ is {\it a
$\star$-valuation overring of } \rm $D$ \rm provided $F^\star
\subseteq FV$,  for each $F \in \boldsymbol{f}(D)$.  Note 
that,
by definition,  the $\star$-valuation overrings coincide with the
$\stf$-valuation overrings. Recall that  \it the $\star$-closure of \rm  
$D$,  defined by:
$$ 
D{\cl} := \bigcup \{(F^\star: F^\star) \mid F \in 
\boldsymbol{f}(D) \}$$
 is an integrally closed 
overring 
of $D$ and, more precisely,  $D\cl = \bigcap \{ V \mid  V$ is a $\star$-valuation 
overring of $D\}$.  
Finally, recall that a valuation overring $V$ of
$D$ is a $\stt$-valuation overring of $D$ if and
only if $V$ is an overring of $D_P$,  for some $\,P\in \QMax^{\stf}(D)$.
For more details on this subject  and for the proofs of the results recalled above,   see   \cite{om},   \cite{HK1}, \cite{HK2}, \cite[Proposition 3.2 and Corollary 3.6]{FL2} and \cite[Theorem 3.9]{fl}.

\begin{corollary}  Let $\star$ be a semistar operation on an integral domain $D$ and let $V$ be a valuation overring of $D$. 
The following statements are equivalent:

\begin{enumerate}

\item[(i)]  $V$ is a $t$-linked overring to $(D, \star)$.

\item[(ii)] $V = V^{\stt}$.

\item[(iii)] $V $ is a $\stt$-valuation overring to $D$.

\end{enumerate}

\end{corollary}
\begin{proof} Note that the $t$-operation on $V$ coincides with the $d$-operation and so (ii)$\Leftrightarrow$(iii) by \cite[Lemma 2.7]{eBFP}.  The equivalence (i)$\Leftrightarrow$(ii) is a particular case of Lemma \ref{t-*-linked} (i)$\Leftrightarrow$(v).  \end{proof}

\begin{remark}  \rm In relation with the previous corollary note that,  given a semistar operation $\star$ on an integral domain $D$, it is known that each overring $R$ of $D$ is $t$-linked to $(D, \star)$ if and only if each valuation overring $V$ of $D$ is $t$-linked to $(D, \star)$
(cf. \cite[Theorem 3.9]{eBF} and \cite[Theorem 2.15]{pt}).

\end{remark}

\begin{lemma}\label{starinc}
 Let $\star$ be a semistar operation on an integral domain $D$.   Then the following
statements are equivalent.
\begin{enumerate}
 \item[(i)]   $D \subseteq K$ is  a  $\stf$-primitive extension (or, a $\star$-primitive extension).  
 \item[(ii)]   $D$ is a  $\stf$-INC-domain.  
 \item[(iii)]    Each  $t$-linked overring to $(D,\star)$   is a  $\stf$-INC  extension of
$D$.
 \item[(iv)]   For each    quasi-$\stf$-prime (or quasi-$\stf$-maximal)  ideal $P$ of $D$,
$D_P \subseteq K$ is a primitive extension.
 \item[(v)]     For each   quasi-$\stf$-prime (or quasi-$\stf$-maximal)  ideal $P$ of $D$,
$D_P$ is an INC-domain.
 \item[(vi)]      For each   quasi-$\stf$-prime (or quasi-$\stf$-maximal)  ideal $P$ of $D$,
$D_P$ is a quasi-Pr\"ufer domain.
\end{enumerate}
\end{lemma}

\begin{proof}
(i)$\Rightarrow$(iv)   Let $P$ be a nonzero quasi-$\stf$-prime  ideal of $D$.  By assumption, if $0 \neq u \in K$, then there is a
polynomial $0 \neq  g \in D[X]$ such that  $\co_D(g)^\star = D^\star$ and
$g(u)=0$. Clearly $g \in D_P[X]$ and $\co_D(g) \nsubseteq P$.  So
$\co_{D_P}(g) = \co_D(g)D_P = D_P$,   and thus $u$ is primitive over $D_P$.

 (iv)$\Rightarrow$(i)  Let $0 \neq u \in K$, and let $I$ be the  (nonzero) 
ideal of $D$ generated by the polynomials $f \in D[X]$ such that
$f(u) = 0$. If   $\co_D(I)^{\stf} = D^{\star}$,   
 there are nonzero polynomials $f_1, f_2,
\dots , f_k \in D[X]$ such that $f_i(u)=0$,  for each $i$,  and    $(\co_D({f_1}), \co_D(f_2), \dots,
\co_D(f_k))^{\star} = D^\star$.    Let   $g := f_1 + X^{n_1+1}f_2 + X^{n_1+n_2+2}f_3 +
\dots + X^{n_1+n_2 + \dots + n_{k-1} + k-1}f_k$,   where $n_i$ is
the degree of $f_i$. Then,  clearly,    $g(u) = 0$ and   $\co_D(g) =(\co_D({f_1}), \co_D(f_2), \dots,
\co_D(f_k))$, thus $\co_D(g)^\star = D^\star$.      So $u$ is $\star$-primitive over $D$. 
 In order to conclude,   it
remains  to show that  $\co(I)^{\stf}=D^\star$. 
 \  Assume that,  for each  $P \in 
\QMax^{\stf}(D)$,  $D_P$ is a   primitive   extension, thus there is a polynomial 
$0 \neq h \in D_P[X]$ such that $h(u)=0$ and $ \co_{D_{P}}(h)  = D_P$. Let
$0 \neq s \in D\setminus P$ with $sh \in D[X]$. Then $ \co_D(sh) 
\nsubseteq P$   (otherwise $ D_P =sD_P = s\co_{D_P}(h) = \co_{D_P}(sh) = \co_D(sh)D_P  \subseteq PD_P$, a contradiction).  Clearly $sh \in I$ and  so $\co_D(I)
\nsubseteq P$ for all   $P \in\QMax^{\stf}(D)  $.    Therefore,   $\co_D(I)^{\stf} = D^{\star}$.

 The equivalences  (iv)$\Leftrightarrow$(v)$\Leftrightarrow$(vi)   follow from Theorem \ref{qp} ((1)$\Leftrightarrow$(4)$\Leftrightarrow$(5)).  

(ii)$\Rightarrow$(iii)  is obvious. 

 (iii)$\Rightarrow$(v)  Let $P$ be a quasi-$\stf$-prime of $D$, let $R$ be  an overring of $D_P$, and let 
$Q_1$ and $Q_2$ are prime ideals of $R$
such that ${Q}_1 \cap D_P = {Q}_2 \cap D_P$, we want to show that  ${Q}_1$ and ${Q}_2$ are
incomparable.  
Let $I$
be  a nonzero finitely generated ideal of $D$ with $I^\star =  D^\star$. Note that 
$I \nsubseteq P$, since $P$ is a quasi-$\stf$-ideal, and hence $D_P = ID_P
\subseteq IR \subseteq R$,  and so   $IR = R$. Thus   $(IR)^{t_R} = R$  and hence  $R$
is $t$-linked to $(D, \star)$.
 By assumption, $R$ is a $\stf$-INC extension of $D$ and  ${Q}_1 \cap D =  {Q}_2 \cap D \subseteq P$, with $P^{\stf} \subsetneq D^{\stf}= D^\star$,  hence ${Q}_1$ and ${Q}_2$ are
incomparable.

 (v)$\Rightarrow$(ii)   Let $R$ be an overring of $D$, and let $Q_1
\subsetneq Q_2$ be prime ideals of $R$ such that
   $Q_2 \cap D $ is contained a quasi-$\stf$-prime $P$ of $D$.   We want to show that $Q_1 \cap D \subsetneq Q_2 \cap D$.  If we consider the extension   $D_P \hookrightarrow
R_{D \setminus P}$  we have  $Q_1R_{D \setminus P} \subsetneq Q_2R_{D
\setminus P} \subsetneq R_{D \setminus P}$ and $D_P$ is an
INC-domain,   by assumption.  Hence $Q_1R_{D \setminus P} \cap D_P \subsetneq
Q_2R_{D \setminus P} \cap D_P$, and thus $Q_1 \cap D = Q_1R_{D
\setminus P} \cap D \subsetneq Q_2R_{D \setminus P} \cap D = Q_2
\cap D$.
\end{proof}

In Theorem \ref{qp} we gave several characterizations of quasi-Pr\"ufer domains. The main goal of this section is to give  a semistar analog characterization theorem for $\stf$-quasi-Pr\"ufer domains, completing the work initiated in Lemma  \ref{*qp}. 
We start with a lemma that extends  to the semistar integral closure  the semistar  operation   versions of the   Cohen-Seidenberg   properties  GU, INC and LO \cite[page 28]{kap}.   (See \cite[Corollary 4.2]{chang2}, \cite[Corollary 1.4]{cz}, or \cite[Theorem 3.3]{wang3} for the star operation versions.)

\begin{lemma} \label{*inc}    Let $\star$ be a semistar operation on an integral domain  $D$ with quotient field $K$.  Let $\overline{D}$ be the integral closure of $D$ (in $K$).  Set $\widetilde{D} := (\overline{D})^{\stt}$, where $\stt$ is the stable semistar of finite type of $D$  associated to $\star$,  and let  $\widetilde{\iota}: D \hookrightarrow \widetilde{D}$ be the canonical embedding.  Set $\ast := (\stt)_{\widetilde{\iota}}$. 

\begin{enumerate}
\item[(a)] $\widetilde{D}$ coincides with the   $\stt$-closure of $D$ (i.e., $\widetilde{D} = \bigcup \{(F^{\stt} : F^{\stt}) \mid F \in \boldsymbol{f}(D) \}$).

\item[(b)]  The inclusion $\widetilde{\iota}: D \hookrightarrow \widetilde{D}$ verifies the properties $\stt$-INC, $\stt$-LO (i.e., for each quasi-$\stt$-prime ideal $P$ of $D$ there exists a $\ast$-prime of $ \widetilde{D}$ such that $Q\cap D = P$), and $\stt$-GU (i.e., if $P \subseteq P'$ are quasi-$\stt$-prime ideals of $D$ and if $Q$ is a $\ast$-prime of $\widetilde{D}$ such that $Q\cap D = P$, then there exists a $\ast$-prime $Q'$ of    $\widetilde{D}$ such that $Q' \cap D = P'$ and $Q \subseteq Q'$).
\end{enumerate}  
\end{lemma}

\begin{proof}   (a)  It is known from \cite[Example 2.1 (c.2)]{fjs} and \cite[Proposition 4.3]{fl} that 
$(\overline{D})^{\stt} =  D{\clt}    = \bigcup \{(F^{\stt} :F^{\stt}) \mid  F \in \boldsymbol{f}(D) \}$, which is an integrally closed overring of $\overline{D}$ (and $D$). 

(b)     Let $P$ be a quasi-$\stt$-prime ideal of $D$.  Consider the prime ideal $PD[X]_{\mathcal{N}^{\stt}}$ and the integral extension $D[X]_{\mathcal{N}^{\stt}} \hookrightarrow \overline{D}[X]_{\mathcal{N}^{\stt}}$. By lying-over, we can find a prime  ideal $\boldsymbol{\mathfrak{Q}}$ in $\overline{D}[X]_{\mathcal{N}^{\stt}}$  such that 
   $\boldsymbol{\mathfrak{Q}} \cap   D[X]_{\mathcal{N}^{\stt}} = PD[X]_{\mathcal{N}^{\stt}}$. Set $Q :=  \boldsymbol{\mathfrak{Q}}  \cap  \widetilde{D} \subseteq \overline{D}[X]_{\mathcal{N}^{\stt}} \cap K = ( \overline{D})^{\stt} = \widetilde{D} $.
   It is easy to see that    $Q$ is a prime ideal of $\widetilde{D}$ such that   $Q^{\ast} = Q^{\stt} = Q$ and $Q \cap D = P$.
   \  Similar arguments prove that $\widetilde{\iota}: D \hookrightarrow \widetilde{D}$ verifies $\stt$-INC and $\stt$-GU.  
\end{proof}

 A domain $D$ is called \it a  Pr\"ufer $\star$--multiplication
domain \rm (for short, \it P$\star$MD\rm) if each nonzero finitely
generated ideal is $\star_{_{\! f}}$-invertible  (cf. for instance
 \cite{fjs} and, for the case of the star operations, \cite{hmm}).  When $\star=v$ we have the
classical notion of P$v$MD   (cf. for instance \cite{griffin},
\cite{mz} and \cite{kang});   when $\star =d$, where $d$ is  the identity (semi)star operation,  we have the notion of
Pr\"ufer domain   \cite[Theorem 22.1]{gilmer}. It is obvious that the notions of P$\star$MD and P$\stf$MD  coincide and it is known that they also coincide with the notion of P$\stt$MD  \cite[Proposition 3.3]{fjs}.  Moreover, when $\star$ is a (semi)star operation then $D$ is a P$\star$MD if and only if $D$ is a P$v$MD and $\stt = t $ (and so $\stt = \stf = t = w$) \cite[Proposition 3.4]{fjs}.   Examples of P$v$MDs that are not  P$\star$MDs (for some (semi)star operation $\star$ on $D$) are given in \cite[Example 3.4]{fjs}.

\begin{theorem}\label{um*}
Let $\star$ be a semistar operation on an integral domain  $D$ with quotient field $K$.  Let $\overline{D}$ be the integral closure of $D$ (in $K$)
  Then the following statements are equivalent.
\begin{enumerate}

 \item[$(1_{\stf})$] $D$ is a  $\stf$-quasi-Pr\"ufer domain.
  
\item[$(4_{\stf})$]  $ D \subseteq K$ is a   $\stf$-primitive   extension.

\item[$(5_\stf)$] $D$ is a  $\stf$-INC-domain.

\item[$(6_{\stt})$]  Set $\widetilde{D}= (\overline{D})^{\stt}$ and let  $\widetilde{\iota}: D \hookrightarrow \widetilde{D}$ be the canonical embedding, then $\widetilde{D}$ is a P${(\stt)}_{\widetilde{\iota}}$MD.  

\item[$(7_{\stf})$]  Each  overring $R$    of $D$ is a   $(\stf)_{\iota}$-quasi-Pr\"ufer domain,   where $\iota: D \hookrightarrow R$ is the canonical embedding.

\item[$(8_{\stf})$] Every prime ideal of $\Na(D, \star_f)$   is extended from $D$.

 \item[$(9_{\stf})$]   $\Na(D, \stf)$    is a quasi-Pr\"ufer domain.
 
 \item[$(10_{\stf})$] The integral closure of   $\Na(D, \stf)$   is a Pr\"ufer domain.
 
 \item[$(11_{\stf})$] $D_P$ is a quasi-Pr\"ufer domain, for each \  quasi-$\stf$-maximal ideal (or,  quasi-$\stf$-prime ideal)  $P$  of $D$.
 
\end{enumerate}

 Moreover, if we assume that $\star$ is a (semi)star operation on $D$, then the previous conditions are also equivalent to the following:

\begin{enumerate}

 \item[$(12_{\stf})$]  Each $t$-linked overring  to   $(D, \stf)$ is a   $t$-quasi-Pr\"ufer domain (or a UM\ \!$t$-domain)   and each $\stf$-maximal ideal   of $D$ is a $t$-ideal.
 
\item[$(13_{\stf})$]  $D$ is  a   $t$-quasi-Pr\"ufer domain (or a UM\ \!$t$-domain)   and each  $\stf$-maximal ideal  of $D$ is a $t$-ideal.

 \item[$(14_{\stf})$] $D$ is a   $t$-quasi-Pr\"ufer domain (or a UM\ \!$t$-domain)  and $\widetilde{\ \stf} = w$. 

\end{enumerate}
\end{theorem}

\begin{proof}

$(1_{\stf})\Rightarrow(4_{\stf})$ Let $0 \neq u \in K$, and let $\ell := X-u$.
Then $Q_\ell := \ell K[X] \cap D[X]$ is a prime ideal of $D[X]$  (since $\ell  \in K[X]$ is irreducible)  and $Q_\ell$  is an upper to zero.  So,  by assumption,  there is a $0 \neq g \in Q_\ell$ such that
$\co_D(g)^\star = D^\star$.    Note that    $g \in Q_\ell \subseteq \ell K[X]$,  so $g = \ell h$ 
for some $h \in K[X]$.  Hence $g(u) = (\ell h)(u) = \ell(u)h(u) = 0$, and
thus $u$ is  $\star$-primitive over $D$.

The equivalences $(4_{\stf})\!\Leftrightarrow\!(5_{\stf})\!\Leftrightarrow\!(11_{\stf})$  are proven in Lemma \ref{starinc}   ((i)$\Leftrightarrow$(ii)$\Leftrightarrow$(vi)).

  $(11_{\stf})\Rightarrow(1_{\stf})$   Let $ \ell  \in D[X]$ be a nonzero polynomial of $D[X]$,  irreducible in $K[X]$,  and  let $Q_\ell   := \ell K[X] \cap D[X]$. Note that  $Q_\ell$ is a prime ideal
of $D[X]$,    $Q_\ell$ is an upper to zero and all upper to zero in $D[X]$ are of this form  \cite[Theorem 36]{kap}.    It is easy to see that,   for each
  quasi-$\stf$-maximal ideal $P$ of $D$, the ideal $ Q_{\ell,P} : =    \ell K[X] \cap D_P[X]$ is a prime ideal
of $D_P[X]$ such that $Q_{\ell,P}  \cap D_P = (0)$ and $Q_{\ell,P}  \cap D[X] =
Q_{\ell} $.   Since $D_P$ is quasi-Pr\"ufer,  $Q_{\ell,P}$ contains a polynomial $h \in D_P[X]$ such that $\co_{D_P}(h) 
\nsubseteq PD_P$  (Theorem \ref{qp} (1)$\Rightarrow$(2)).   Choose $s \in D \setminus P$ with $sh \in D[X]$.  Note that  $sh  \in Q_{\ell,P} \cap D[X] = Q_\ell$
and that 
  $s\co_D(h) = \co_D(sh) \nsubseteq P$,  because $\co_{D_P}(sh) 
\nsubseteq PD_P$.  Since the last property holds for each  quasi-$\stf$-maximal ideal $P$ of $D$,  then $\co_D(Q_\ell)^{\stf} = D^\star$.  We conclude that $D$ is a $\stf$-quasi-Pr\"ufer domain by   Lemma \ref{*fqp}.   

 $(1_{\stf})\Rightarrow(8_{\stf})$  Suppose that $D$ is a   $\stf$-quasi-Pr\"ufer domain,   and let
 $\boldsymbol{\mathfrak{Q}}$   be a prime ideal of  $\Na(D, \star)$.  Then there is a prime ideal
$Q$ of $D[X]$ such that  $\boldsymbol{\mathfrak{Q}} = Q\Na(D, \star)=QD[X]_{\mathcal{N}^\star}$  and so   $Q \cap \mathcal{N}^\star =
\emptyset$. Let $P := Q \cap D$.     If $P[X] \subsetneq Q$, pick ~$ q
\in Q \setminus P[X]$, and let   $Q_1$  be an upper to zero in $D[X]$
such that ~$q \in Q_1\subseteq Q$   \cite[Theorem A]{dl}.   Since $D$
is a  $\stf$-quasi-Pr\"ufer domain   and  $Q_1$   is an upper to zero,  there is a nonzero polynomial     $g \in Q_1$  such that $\co_D(g)^\star = D^\star$,
and hence $g \in \mathcal{N}^\star \cap Q$,  a contradiction. So $Q = P[X]$, and
thus $\boldsymbol{\mathfrak{Q}}  =  P\Na(D, \star)$.

$(8_{\stf})\Rightarrow(1_{\stf})$ Suppose that $D$ is not a  $\stf$-quasi-Pr\"ufer domain. By Lemma  \ref{*fqp} ($(1_{\stf})\Leftrightarrow(2_{\stf})$)    then
there is an upper to zero $Q$ in $D[X]$ such that  $Q \cap \mathcal{N}^\star=
\emptyset$. Hence $QD[X]_{\mathcal{N}^\star} =Q\Na(D, \star)$ is a proper prime ideal of
$\Na(D, \star)$.   Note that  $ Q\Na(D,\star) \neq P\Na(D,\star)$ for all nonzero  prime ideals
$P$ of $D$,  since $Q$ is an upper to zero. This fact contradicts the assumption $(8_{\stf})$.

$(9_{\stf})\Leftrightarrow(10_{\stf})$  follows from    Theorem  \ref{qp} ((1)$\Leftrightarrow$(6)).

$(9_{\stf})\Rightarrow(11_{\stf})$  Let $P$ be a  quasi-$\stf$-maximal ideal  of $D$. Then
$P\Na(D, \star)$   is a maximal ideal of   $\Na(D, \star)$  \cite[Proposition 1.3 (3)]{fl}. 
Hence  $\Na(D_P) = D_P[X]_{PD_P[X]} $ $= (D[X]_{\mathcal{N}^\star})_{PD[X]_{\mathcal{N}^\star}} =\Na(D, \star)_{P{\rm \footnotesize Na}(D, \star)}$  (cf. also \cite[Theorem 3.8]{fl}).  Since we are assuming that $\Na(D, \star)$ is quasi-Pr\"ufer, then $\Na(D, \star)_{P\mbox{\rm \footnotesize Na}(D, \star)}= \Na(D_P)$ is quasi-Pr\"ufer  
and thus $D_P$ is a quasi-Pr\"ufer domain   (by 
Theorem \ref{qp} ((9)$\Rightarrow$(1)). 

$(11_{\stf})\Rightarrow(9_{\stf})$  Let $Q \in  \Max(\Na(D, \star))$. We know that $Q = P\Na(D, \star)$, for some $P \in \QMax^{\stf}(D)$  and that  $\Na(D, \star)_Q = \Na(D_P)$  \cite[Proposition 1.3]{fl}. Therefore,   if $D_P$ is quasi-Pr\"ufer, then $\Na(D, \star)_Q = \Na(D_P)$ is quasi-Pr\"ufer (Theorem  \ref{qp} (1)$\Rightarrow$(9)). Thus $\Na(D, \star)$ is quasi-Pr\"ufer (Theorem \ref{qp} (11)$\Rightarrow$(1)).   
 
$(1_{\stf})\Rightarrow(7_{\stf})$   Let $R$ be an overring to $(D, \stf )$ and,  for simplicity of notation,  set $\ast := (\stf\! )_\iota$.   The property $(7_{\stf})$  holds if we show that $R_Q$ is a quasi-Pr\"ufer domain for all $Q \in \QMax^{\ast}(R)$,   since we already proved  that $(1_{\stf})\Leftrightarrow(11_{\stf})$.      Note that the prime ideal $P:=Q \cap D$  is such that $P \subseteq P^{\stf} \cap D \subseteq Q^{\stf} \cap D = 
 (Q^{\stf} \cap R) \cap D = (Q^{\ast} \cap R) \cap D =Q \cap D =P$. Since $P$ is a quasi-$\stf$-prime ideal of $D$ and, by assumption, $D$ is a  $\stf$-quasi-Pr\"ufer domain,   then $D_P$ is quasi-Pr\"ufer (by $(1_{\stf})\Rightarrow(11_{\stf})$).  Therefore $R_Q$ , which is an overring of $D_P$, is also quasi-Pr\"ufer (Theorem \ref{qp} (1)$\Rightarrow$(7)). 
 
$(7_{\stf})\Rightarrow(1_{\stf})$ is trivial.
    
$(1_{\stf})\Rightarrow(6_{\stt})$ 
    We already proved that $(1_{\stf})$ is equivalent to $(9_{\stf})$.  Therefore we can assume that   $\Na(D, \stf) \ $ $ (= \Na(D, \stt) $ by \cite[Corollary 3.7]{fl}) is a quasi-Pr\"ufer domain, i.e., the integral closure  $\overline{\Na(D, \stt)}$ of $\Na(D, \stt)$  in $K(X)$ is a Pr\"ufer domain.   Note that $\overline{\Na(D, \stt)} = \overline{D}[X]_{\mathcal{N}^{\stt}} $,
  where  $\mathcal{N}^{\stt}= \mathcal{N}^{\stf} = \{g \in D[X] \mid   g \neq 0 \mbox{ and }   \co_D(g)^{\stt} =D^{\stt} \}$.
 For the sake of simplicity, set $\ast :=  (\stt )_{\widetilde{\iota}}$.   Clearly $\ast$ is a stable (semi)star operation of finite type on $\widetilde{D}$. Moreover  $\Na(\widetilde{D}, \ast)=\widetilde{D}[X]_{\widetilde{\mathcal{N}}}$,  where $\widetilde{\mathcal{N}} := \mathcal{N}^{\ast} = 
   \{h \in \widetilde{D}[X] \mid   h \neq 0 \mbox{ and }  \co_{\widetilde{D}}(h)^\ast = \widetilde{D} \}$.  Since it is clear that $\mathcal{N}^{\stt}$ is also a multiplicative set in $\widetilde{D}[X]$ and that  $\mathcal{N}^{\stt} \subseteq \widetilde{\mathcal{N}} $, then $\overline{\Na(D, \stt)} \subseteq \Na(\widetilde{D}, \ast)$ and so $\Na(\widetilde{D}, \ast)$ is a Pr\"ufer domain   \cite[\ec~Theorem 26.1]{gilmer}.  By \cite[Theorem 3.1 (i)$\Leftrightarrow$(iii)]{fjs},  this  is equivalent to say that $\widetilde{D}$ is a P$\ast$MD.
   
   $(6_{\stt})\Rightarrow(10_{\stf})$ With the notation used in the proof of $(1_{\stf})\Rightarrow(6_{\stt})$,     the present hypothesis is equivalent  to assume that  $\Na(\widetilde{D}, \ast)$ is a  Pr\"ufer domain.   The conclusion will trivially follow if we show that     $\overline{\Na(D, \stt)} = \Na(\widetilde{D}, \ast)$,  i.e., $ \overline{D}[X]_{\mathcal{N}^{\stt}} =\widetilde{D}[X]_{\widetilde{\mathcal{N}}}$\  \!.  
 
  Note that  $ \mathcal{N}^{\stt} = D[X]\setminus \bigcup \{ P[X]  \mid P \in \QMax^{\stt}(D)\}$,
 $ \widetilde{\mathcal{N}} =\widetilde{D}[X] \setminus \bigcup \{ Q[X]   \mid Q \in \Max^{\ast}(\widetilde{D}) \}$,   and  $ \overline{D}[X]_{\mathcal{N}^{\stt}} \subseteq \widetilde{D}[X]_{\widetilde{\mathcal{N}}}$\ \!.  By Lemma \ref{*inc} (b)  the natural embedding \ec~$\widetilde{\iota}: D \hookrightarrow \widetilde{D}$ verifies $\stt$-LO, $\stt$-INC and $\stt$-GU.  It is not difficult to see that  a prime ideal $Q$  of $\widetilde{D}$ belongs to   $ \Max^{\ast}(\widetilde{D})$ if and only if $Q\cap D$  belongs to $ \QMax^{\stt}(D)$.
 
        As a matter of fact, let $Q$ be a prime ideal in $\widetilde{D}$. Assume that $P:= Q\cap D \in \QMax^{\stt}(D)$. By $\stt$-LO we can assume that $Q$ is a $\ast$-prime in $\widetilde{D}$.   Let  $M\in \Max^{\ast}(\widetilde{D})$, such that $M \supsetneq Q$. By $\stt$-INC we have $M\cap D \supsetneq Q \cap D=P$. Therefore $M \cap D \subseteq (M\cap D)^{\stt} \cap D =(M^{\stt} \cap{D}^{\stt}) \cap D \subseteq (M^{\stt} \cap\widetilde{D}) \cap D = (M^{\ast} \cap\widetilde{D}) \cap D = M \cap D$ and so we reach a contradiction (i.e., $P$ is not in $\QMax^{\stt}(D)$). 
     Conversely, let $Q \in \Max^{\ast}(\widetilde{D})$ and assume that $P:= Q \cap D \subsetneq P' $, for some prime ideal $P'$ of $D$  such that  $ P' ={P'}^{\stt} \cap D \in \QMax^{\stt}(D)$.  By $\stt$-GU, there exists a $\ast$-prime ideal $Q'$ of $\widetilde{D}$ such that $Q' \cap D = Q$ and $Q \subsetneq Q'$ and this is a contradiction.
     
    From the fact that   a prime ideal $Q$  of $\widetilde{D}$ belongs to   $ \Max^{\ast}(\widetilde{D})$ if and only if $Q\cap D$  belongs to $ \QMax^{\stt}(D)$, we deduce that the ideals of $ \overline{D}[X]$ that are maximal with the property of being disjoints from ${\mathcal{N}^{\stt}} $ are the ideals $\{ (Q\cap\overline{D})[X] \mid  Q \in \Max^{\ast}(\widetilde{D}) \}$. From this fact it follows easily that
     $ \overline{D}[X]_{\mathcal{N}^{\stt}} =\widetilde{D}[X]_{\widetilde{\mathcal{N}}}$
     
 $(13_{\stf})\Leftrightarrow(14_{\stf})$   The second part of condition $(13_{\stf})$ implies that $\Max^{\stf}(D) = \Max^t(D)$ and so $\stt = w$. Conversely, if $\stt = w$, then 
 $\Max^{\stf}(D) = \Max^{\stt}(D) = \Max^w(D)= \Max^t(D)$, and so each quasi-$\stf$-ideal maximal  of $D$ is a $t$-ideal. 
    
   $(1_{\stf})\Rightarrow(13_{\stf})$  Under the present assumptions   $\stf $ is a (semi)star operation of finite type   on  $D$, then $\stf \leq
t$   (essentially by \cite[Theorem 34.1 (4)]{gilmer}). Therefore   if $D$ is a  $\stf$-quasi-Pr\"ufer domain   then $D$ is also   a   $t$-quasi-Pr\"ufer domain (Corollary \ref{<} (a)).   Let $P$ be a  
$\stf$-maximal ideal  of $D$. Since we already proved that  $(1_{\stf})\Rightarrow(11_{\stf})$,
  then $D_P$ is a quasi-Pr\"ufer domain.    By Corollary \ref{qp-umt}  ((1)$\Rightarrow$(13))  $PD_P$ is a $t$-ideal  in $D_P$   and thus $P = PD_P \cap D$ is a $t$-ideal of $D$ \cite[Lemma
3.17]{kang}.

 $(13_{\stf})\Rightarrow(10_{\stf})$  Since each   $\stf$-maximal ideal of $D$ is a
$t$-ideal, then necessarily $\Max^{\stf}(D) =\Max^t(D)$ and hence   $\mathcal{N}^\star = \mathcal{N}^{\stf} = \mathcal{N}^t = \mathcal{N}^v$. Thus $\Na(D, \star) =\Na(D, v)$ and   $\Na(D, v)$   has Pr\"ufer integral closure by
 \cite[Theorem 2.5]{fgh},  since $D$ is an UM$t$-domain   ($= t$-quasi-Pr\"ufer domain).  

   $(1_{\stf})\Rightarrow(12_{\stf})$   Let $R$ be a $t$-linked overring to $(D, \stf)$, then $R = R^{\stt}$ (Lemma \ref{t-*-linked} ((i)$\Rightarrow$(v))). Let $\iota : D \hookrightarrow R$ be the canonical embedding, then   $(\stt)_\iota$ is a (semi)star  operation of finite type on $R$.   Since $D$ is   a  $\stf$-quasi-Pr\"ufer domain  or, equivalently, a  $\stt$-quasi-Pr\"ufer domain (Corollary \ref{<} (c))  then, by $(1_{\stf})\Leftrightarrow(7_{\stf})$, $R$ is a   $(\stt)_\iota$-quasi-Pr\"ufer domain.    Since in the present situation $(\stt)_\iota \leq t_R $ then  $R$ is also a  $t_R$-quasi-Pr\"ufer domain.  Moreover, we already proved that  $(1_{\stf})\Rightarrow(13_{\stf})$ thus each  $\stf$-maximal ideal  of $D$ is a $t$-ideal.  
      
$(12_{\stf})\Rightarrow(13_{\stf})$  is trivial. \end{proof}

Let $\star$ be a semistar operation on an integral domain $D$. Recall that  a P$\star$MD $D$ can be characterized by the fact that  $D_P$ is a valuation domain for each $P \in \QMax^{\stf}(D)$ \cite[Theorem 3.1]{fjs}. Thus, since a valuation domain is trivially quasi-Pr\"ufer, a P$\star$MD is a $\stf$-quasi-Pr\"ufer domain by Theorem  \ref{um*} ($(1_{\stf}) \Leftrightarrow (11_{\stf}))$.  This fact generalizes the well known property that a P$v$MD is a UM$t$ domain (Corollary \ref{<}  (b)).   However, a P$\star$MD need not be integrally closed (cf. \cite[Example 3.10]{fjs}), while being a P$v$MD is equivalent  to being an integrally closed UM$t$ domain \cite[Proposition 3.2]{hz}. The next corollary gives an appropriate generalization of  the previous result to the case of semistar operations. 

\begin{corollary} \label{p*md-um*}
 {\em   \cite[Theorem 3.2]{fjs} }  Let $\star$ be a semistar operation on an integral domain $D$ with quotient field $K$. Then the following statements are equivalent.
\begin{enumerate}
\item[(i)] $D$ is a P$\star$MD.
\item[(ii)]  $D$ is a $\stf$-quasi-Pr\"ufer domain and $D_P$ is integrally closed for all $P \in \QMax^{\stf}(D)$.
\item[(iii)]  $D$ is a $\stf$-quasi-Pr\"ufer domain and $D^{\stt}$ is integrally closed.
\end{enumerate}
\end{corollary}

\begin{proof}
The implication (i)$\Rightarrow$(ii) was already proved just before the statement of Corollary \ref{p*md-um*}.

(ii)$\Rightarrow$(iii) This follows from \cite[Theorem 52]{kap} because $D^{\stt} = \bigcap \{D_P \mid  P \in \QMax^{\stf}(D)\}$ and each $D_P$   is integrally closed, by assumption.

 (iii)$\Rightarrow$(i) Let $\iota:D \hookrightarrow D^{\stt}$  be the canonical embedding and set $\ast := (\stt)_\iota$ (thus $E^\ast =   E^{\stt}$ for all $E \in \overline{\boldsymbol F}(D^{\stt}) \ (\subseteq \overline{\boldsymbol F}(D)$). 
 Then $\Na(D, \stf) = \Na(D, \stt) = \Na(D^{\stt}, \ast)$ \cite[Corollary 3.5]{fl}.  On the other hand $\Na(D^{\stt}, \ast)$ is integrally closed, because $D^{\stt}$ is integrally closed by assumption, and 
  $\Na(D, \stf)$ is 
  quasi-Pr\"ufer domain by Theorem 2.17 ($(1_{\stf}) \Leftrightarrow (9_{\stf})$). Putting these two pieces of information together we deduce that  $\Na(D, \stf)$ is  a Pr\"ufer domain and thus $D$ is a P$\stf$MD (or, a P$\star$MD) by  \cite[Theorem 3.1]{fjs}. 
\end{proof}

The following corollary follows immediately from Theorem \ref{um*} ($(1_{\stf})\!\Leftrightarrow\!(6_{\stt})$) and from \cite[Proposition 3.4]{fjs}.

     \begin{corollary} \label{umt-closure}  {\em (cf. \cite[Theorem 4.2]{wang3})}  Let $D$ be an integral domain  with quotient field $K$.   Set $\widetilde{D}= (\overline{D})^{w_D}$ and let  $\widetilde{\iota}: D \hookrightarrow \widetilde{D}$ be the canonical embedding.  The following statements are equivalent:
     \begin{enumerate}
     \item [(i)]  $D$ is a UM$t_D$-domain.
      \item [(ii)]  $\widetilde{D}$ is a P$(w_D)_{\widetilde{\iota}}$MD.
   \item[(iii)] $\widetilde{D}$ is a P$v_{\widetilde{D}}$MD and  $(w_D)_{\widetilde{\iota}} =w_{\widetilde{D}}=t_{\widetilde{D}} $. \hfill $\Box$
   \end{enumerate}
     
     \end{corollary}

     \medskip
  
  We have already mentioned in Remark \ref{umt} (d) the  interesting open problem  of establishing whether the integral closure
of a UM$t$-domain is a P$v$MD.  For a negative answer to this problem we need examples of integral domains $D$ such that the integral closure $\overline{D}$ is not $t$-linked to $D$ (Remark  \ref{umt} (a)).  This is not an easy task, even in a general situation.  Note that the integral closure $\overline{D}$ is $t$-linked to $D$ if $D$ is one-dimensional    \cite[Corollary 2.7]{dhlz2}   or if $D$ is quasi-coherent (e.g., $D$ is Noetherian)    \cite[Corollary 2.14 (a)]{dhlz2}.  A first class of examples of   integral domains of dimension $ \geq 3$ such that the integral closure $\overline{D}$ is not $t$-linked to $D$ was given in \cite[Example 4.1]{dhlrz}.    The 2-dimensional case was left open in that paper.   A first example in dimension two was given by Dumitrescu \cite{dumitrescu}, using the $A +XB[X]$ constructions. We give next another example of this type.

\begin{example} \label{example} \sl     A quasi-local strong Mori non Noetherian 2-dimensional UM$t$-domain $D$ such that $\overline{D}$ is not $t$-linked to $D$ but still $\overline{D}$ is   a P$v_{\overline{D}}$MD.  \rm

For this purpose we use a construction due to Heinzer, Ohm and Pendleton \cite[Example 2.10]{hop}.  Let $K$ be a field, $X, Y$ indeterminates over $K$, let $V$ be  the $X$-adic valuation ring of $K(X, Y)$, i.e.  $V := K(Y)[X]_{(X)}$, and let $M_X :=  XK(Y)[X]_{(X)} $ be the maximal ideal of $V$ (hence $V = K(Y)+M_X$). Also, let $D_1 := K[X, Y]_{(X, Y)}$,  $M_1:= (X,Y)K[X, Y]_{(X, Y)}$,   $k_T := K(Y+ \frac{1}{Y}) \subsetneq K(Y)$ and set $T := k_T +M_X$, and $D := T \cap D_1$.   Note that if we consider the Krull overring $R := D_1[1/X] = \bigcap \{ {D_1}_{P_1} \mid P_1 \neq (M_X\cap D_1) \mbox { with $P_1$ height 1 prime ideal of $D_1$}\}$ of $D$ (and of $D_1$) \ec~\cite[Corollary 1.5 and Proposition 3.15]{fossum}, then we also have $D = R \cap T$ (and ($D_1 = R \cap V$).

\begin{enumerate}
\item[(a)] $T$ is a 1-dimensional Noetherian pseudo-valuation domain (or, PVD) with maximal ideal $M_X$ and associate valuation overring $V$. Moreover  the integral closure $\overline{T}$ of $T$  coincides with $V$.
\end{enumerate}

Note that $k_T \hookrightarrow K(Y)$ is a finite extension, since $Y$ is a root of the polynomial $ Z^2 -\left((Y^2 +1)/Y\right)Z +1$ in the indeterminate $Z$ with coefficients in $k_T$. The conclusion follows
from \cite[Theorem 3.1 and Corollary 3.4]{hh}.

\begin{enumerate}
\item[(b)]  Let  $Q := M_X \cap D = XK[X, Y]_{(X, Y)} = M_X \cap D_1$. Then $D_1 \subsetneq V= (D_1)_Q$ and that $D$ and $D_1$ have a common  prime ideal, i.e., $Q$. In particular, the map $H_1 \mapsto H := H_1 \cap D$ establishes   a 1--1 correspondence betweeen the prime ideals  of $D_1$ not containing $Q$ and the prime ideal of $D$ not containing $Q$ and, moreover, 
$D_{H} = (D_1)_{H_1}$.  For the remaining localization of $D$ at the prime $Q$, we have   $D_Q = T \subsetneq (D_1)_{Q} $.
\end{enumerate}

After remarking that $Q$ is a common ideal of $D$ and $D_1$, the first part follows from the general properties of the pullback diagrams \cite[Claim (c) in the proof of Theorem 1.4]{f}. The last statement is proved in \cite[Lemma, page 152]{hop}.  

\begin{enumerate}
\item[(c)]  $D$ is a quasi-local domain with maximal ideal $M:= M_1 \cap D$, with complete integral closure equal to $D_1$ and dim$(D)=2$.  
\end{enumerate}

The first part of the statement  is proved in \cite[Example 2.10, page 152]{hop}.  The reamining part follows from the fact that $D$ and $D_1$ have $Q$ as common ideal \cite[Lemma 26.5]{gilmer} and from the fact that $\dim(D_1) =2$.

\begin{enumerate}
\item[(d)]   $D$ is  a strong Mori domain   with dim$^t(D) = 1$,
 \end{enumerate}
 
 Let $\Lambda_1 :=\{ P_1 \in D_1 \mid P_1 \mbox{ is a height 1 prime ideal of } D_1,  \,  P_1 \neq Q \}$ (resp.,   $\Lambda  :=\{ P  \in D \mid P  \mbox{ is a height 1 prime ideal of } D,  \,  P \neq Q \}$). From (b) and from the presentations  $D=  \left(\bigcap \{ {D}_{P} \mid P  \in \Lambda\}\right)  \cap D_Q $
 $=$
 $ \left(\bigcap \{ {D_1}_{P_1} \mid P_1 \in \Lambda_1\}\right)  \cap D_Q $ 
 $=$ 
  $(\bigcap \{ {D_1}_{P_1} \mid$  $ P_1 \in \Lambda_1\})  \cap T$ $\subsetneq$ 
 $\left(\bigcap \{ {D_1}_{P_1} \mid P_1 \in \Lambda_1\}\right)  \cap V$
 $=$ 
 $  \left(\bigcap \{ {D_1}_{P_1} \mid P_1 \in \Lambda_1\}\right)  \cap {D_1}_Q  = D_1$ \ we deduce  that $D$ is a Mori domain  (in particular, $t=v$) \cite[Construction 4.1 and Theorem 4.3]{bg}.  Obviously, all the height 1 prime ideals of $D$ are $t$-ideals of $D$, but  the maximal ideal $M$ is not a $t$-ideal (or a $v$-ideal) of $D$  \cite[Theorem 4.3 (f)]{bg}.   Henceforth dim$^t(D) = 1$ and $\Max^t(D) = \{P\in  \Spec(D) \mid \hgt(P) =1 \}$. Furthermore,  note that $D_P$ is Noetherian for all $P \in \Max^t(D)$  and  each nonzero element of $D$ lies in only finitely many maximal $t$-ideals of $D$ (because this property holds in $D_1$) \cite[Theorem 4.3 (a)]{bg}.  Therefore, by \cite[Theorem 1.9]{WMc99},  $D$ is a strong-Mori domain (i.e., $D$ verifies the acc on the $w$-ideals \cite{WMc97}) and, clearly,  $D = \bigcap \{ D_P \mid P \in \Max^t(D)  \} = D^w$.

\begin{enumerate}
\item[(e)]  $D$ is a UM$t$-domain.
\end{enumerate}
By (d) dim$^t(D) = 1$, then  $D$ is a UM$t$-domain by  \cite[Corollary 3.2 ((6)$\Rightarrow$(1))]{cz}.

\begin{enumerate}
\item[(f)] The integral closure $\overline{D}$ of $D$   coincides with $(W_1 \cap W_2) + Q$, where $W_1 := K[Y]_{(Y)}$ and $W_2 := K[\frac{1}{Y}]_{(\frac{1}{Y})}$. Therefore $\overline{D} \subsetneq D_1$, 
$\overline{D} / Q = W_1 \cap W_2$ is a  semi-quasi-local PID with two maximal ideals and $D_1/Q = W_1$.
\end{enumerate}

The first part of the statement is proved in  \cite[Example 2.10, page 152]{hop}.  The remaining part is an easy consequence of the first part \cite[Theorem 107]{kap}.

The following three statements are immediate consequences of (f).

\begin{enumerate}
\item[(g)]  $\overline{D}$ is semi-quasi-local with two maximal $\overline{M}_1$ and $\overline{M}_2$ such that $\overline{M}_1 \cap D = \overline{M}_2 \cap D = M$ and ht$(\overline{M}_1)$ = ht$(\overline{M}_2) = 2$. Moreover,  $\overline{D}_{\overline{M}_1}/Q\overline{D}_{\overline{M}_1} = W_1$ and 
$\overline{D}_{\overline{M}_2}/Q\overline{D}_{\overline{M}_2} = W_2$.

\item[(h)]  The only prime ideals of $\overline{D}$ containing $X$ (i.e., the prime ideal $Q=XD_1$) are $\overline{M}_1$, $\overline{M}_2$ and, obviously, $Q$. 

\item[(i)]   $D$ and $\overline{D}$ have a common  prime ideal $Q $, then --as in point (b)-- the map $\overline{H}\mapsto H := \overline{H}\cap D$ establishes   a 1--1 correspondence betweeen the  prime ideals  of $\overline{D}$  not containing $Q$ and the prime ideal of $D$ not containing $Q$ and, moreover, 
$D_{H} =  \overline{D}_{\overline{H}}$. Furthermore, as a consequence of (a) and (b), $\overline{D}_Q = V$.

\item[(j)]  $\overline{D} \subsetneq  (\overline{D})^{w_D}   = D_1$. Therefore $\overline{D}$ is not $t$-linked to $D$ (Lemma \ref{t-*-linked} ((i)$\Rightarrow$(v))) and $D$ is not Noetherian.

  \end{enumerate}

As already remarked in \cite[Example 2.10, page 152]{hop},  we have $(\overline{D})^{w_D}  = \bigcap \{\overline{D}  \cdot\!D_P \mid P \in \Max^t(D) \} = 
\overline{D}_{D\setminus Q} \cap (\bigcap \{\overline{D}_{D\setminus P}  \mid P \in  \Lambda\}) = 
V \cap  (\bigcap \{ {D_1}_{ P_1}  \mid P_1 \in \Lambda_1\})= V \cap R  = D_1$ (cf. also \cite[Theorems 1.3 and 3.1]{cz}).  The claim that $D$ is not Noetherian is a consequence of the fact that $\overline{D} \neq D_1$ and that, by (c),  $D_1$ is the complete integral closure of $D$.

\smallskip

 Set  $A:= W_1 \cap W_2$, $B :=W_1$, and let $\mathfrak{m_1}$, $\mathfrak{m_2}$ be the maximal ideals of $A$, with $A_{\mathfrak{m_1}} =W_1$ and $A_{\mathfrak{m_2}} =W_2$ (cf. (f)). By the previous considerations, we have the following   pullback   diagrams of canonical homomorphisms:

$$
\begin{CD}
D  @>>>   D/Q \\
 @VVV       @VVV   \\
 \overline{D}  @>>> \overline{D}/Q  @= A\\
 @VVV       @VVV   @VVV \\
 D_1= (Q:Q) @>>> D_1/Q @=B
\end{CD}
$$

\begin{enumerate}
\item[(k)]  $\overline{D} $ is a P$v_{\overline{D}}$MD.
\end{enumerate}

We claim that for each  prime $t$-ideal   $\mathfrak{p}$ of $A$ either $A_{\mathfrak{p}}$ is a valuation domain and $B_{A\setminus \mathfrak{p}}$ is a field or there exists a finitely generated ideal $\mathfrak{f}$ of $A$, $\mathfrak{f} \subseteq \mathfrak{p}$ such that $(A: \mathfrak{f}) \cap A_{\mathfrak{p}} =A$. 
As a matter of fact, by (f),  $A$ is a PID with $\Max(A) = \{\mathfrak{m_1}, \mathfrak{m_2} \}$, then the set of   prime $t$-ideals  of $A$ coincides with $\Max(A)$.  Clearly, $A_{\mathfrak{m_2}} =W_2$ and $B_{A \setminus \mathfrak{m_2}} = (W_1)_{A \setminus \mathfrak{m_2}}$ is the quotient field of $B$ (and of $A$).
On the other hand $A_{\mathfrak{m_1}} =W_1= B_{A \setminus \mathfrak{m_1}} $, but if $\mathfrak{m_1} = \pi A$, then $(A: \pi A) \cap A_{\pi A} = \pi^{-1}A \cap A_{\pi A} =A$, since 
$A \cap \pi A_{\pi A} =   \pi A$.
 \ Now the statement  follows from \cite[Theorems 4.8 and 5.2]{ht}.

\end{example}

\begin{remark} \label{um*D-bar}    \rm  (a)  With the notation of Theorem  \ref{um*},
let $D\overset{\overline{\iota}}{\hookrightarrow} \overline{D}$, \ $\overline{D}\overset{j}{\hookrightarrow} \widetilde{D}$ and $D\overset{\widetilde{\iota}}{\hookrightarrow} \widetilde{D}$ be the canonical embeddings and so $\widetilde{\iota} = j \circ \overline{\iota}$.     Note that  \it   the statement  $(6_{\stt})$ is equivalent to each of the following: \rm 
\begin{enumerate}
\item[$({6^\prime}_{\stt})$]  \it $\widetilde{D}$ is a P${v}_{\widetilde{D}}$MD and  ${(\stt)}_{\widetilde{\iota}} = w_{\widetilde{D}} = t_{\widetilde{D}}$. 

\item[$(\overline{6}_{\stt})$]  \it $\overline{D}$ is a  P${(\stt)}_{\overline{\iota}}$MD.  

\item[$(\overline{6}^{\prime}_{\stt})$]  \it $\overline{D}$ is a P${({v}_{\widetilde{D}})}^{{j}}$MD and
${(\stt)}_{\widetilde{\iota}} = w_{\widetilde{D}} = t_{\widetilde{D}}$. 
\rm 
\end{enumerate}

\indent  The equivalence  $({6}_{\stt}) \Leftrightarrow ({6^\prime}_{\stt})$ follows immediately from  \cite[Proposition 3.4]{fjs}, since ${(\stt)}_{\widetilde{\iota}} $ is a (semi)star operation on $\widetilde{D}$.

$({6}_{\stt}) \Rightarrow ({\overline{6}}_{\stt})$ Set
$$ 
\begin{array}{rl}
\mathcal{N}^{\stt} =& \hskip -5 pt  \{g \in D[X] \mid  g \neq 0 \mbox{ and }  \co_D(g)^{\stt} =D^{\stt}=\widetilde{D} \}, \\

\overline{\mathcal{N}} :=&  \hskip -5 pt  \mathcal{N}^{{(\stt)_{\overline{\iota}}}}= 
   \{\ell \in \overline{D}[X] \mid    \ell \neq 0 \mbox{ and }  \co_{\overline{D}}(\ell)^{{(\stt)_{\overline{\iota}}}}=  \overline{D}^{{(\stt)_{\overline{\iota}}}} = \widetilde{D} \}, \\
   
\widetilde{\mathcal{N}} :=&  \hskip -5 pt \mathcal{N}^{{(\stt)_{\widetilde{\iota}}}}= 
   \{h \in \widetilde{D}[X] \mid    h  \neq 0 \mbox{ and }  \co_{\widetilde{D}}(h)^{{(\stt)_{\widetilde{\iota}}}}= \widetilde{D} \}.

\end{array}
$$
Clearly  $\mathcal{N}^{\stt} \subseteq \overline{\mathcal{N}}   \subseteq \widetilde{\mathcal{N}} $,  in particular, $ \overline{D}[X]_{\mathcal{N}^{\stt}} \subseteq  \overline{D}[X]_{\overline{\mathcal{N}} } \subseteq \widetilde{D}[X]_{\widetilde{\mathcal{N}} }$\ \!.  On the other hand, $ \Na(\overline{D}, {{(\stt)_{\overline{\iota}}}}) = \overline{D}[X]_{\overline{\mathcal{N}} }$ and 
$\Na(\widetilde{D}, {{(\stt)_{\widetilde{\iota}}}}) = \widetilde{D}[X]_{\widetilde{\mathcal{N}} }$\ \!.
Recall that in the proof  $(6_{\stt})\Rightarrow(10_{\stf})$ of Theorem \ref{um*}, we have shown that  $\overline{D}[X]_{\mathcal{N}^{\stt}} =\widetilde{D}[X]_{\widetilde{\mathcal{N}}}$\ \!.  Therefore, in particular, 
$\Na(\overline{D}, {{(\stt)_{\overline{\iota}}}}) =\Na(\widetilde{D}, {{(\stt)_{\widetilde{\iota}}}})$. Henceforth, if 
$(6_{\stt})$ holds then  $\Na(\overline{D}, {{(\stt)_{\overline{\iota}}}}) \ (=\Na(\widetilde{D}, {{(\stt)_{\widetilde{\iota}}}}))$ is a Pr\"ufer domain and so $\overline{D}$ is a  P${(\stt)}_{\overline{\iota}}$MD   \cite[Theorem 3.1 (i)$\Leftrightarrow$(iii))]{fjs}.

   $({\overline{6}}_{\stt}) \Rightarrow ({6}_{\stt})$  
    By assumption and by \cite[Theorem 3.1 (i)$\Leftrightarrow$(iii))]{fjs}  $\Na(\overline{D}, {{(\stt)_{\overline{\iota}}}})$ is a Pr\"ufer domain  and then obviously each overring of   $\Na(\overline{D}, {{(\stt)_{\overline{\iota}}}})$   is a Pr\"ufer domain. In particular $\Na(\widetilde{D}, {{(\stt)_{\widetilde{\iota}}}})$ is a Pr\"ufer domain and thus  $(6_{\stt})$ holds again by  \cite[Theorem 3.1 (i)$\Leftrightarrow$(iii))]{fjs}.

   $ ({\overline{6}^\prime}_{\stt})  \Leftrightarrow ({\overline{6}}_{\stt}) $ Note that, for each $E \in \overline{\boldsymbol{F}}(\overline{D})$, we have:
$$
E^{(\stt)_{\overline{\iota}}} = E^{\stt} = \bigcap \{ED_P \mid P \in \QMax^{\stt}(D) \}  =  (E\widetilde{D})^{(\stt)_{\widetilde{\iota}}} =  E^{((\stt)_{\widetilde{\iota}})^j}\,.$$ 
Therefore $ {(\stt)_{\overline{\iota}}} ={((\stt)_{\widetilde{\iota}})^j}$.   Henceforth it is straightforward that $ ({\overline{6}^\prime}_{\stt})  \Rightarrow ({\overline{6}}_{\stt}) $ after recalling  that 
$((v_{\widetilde{D}})^j)_{\!_f} =  (t_{\widetilde{D}})^j$.  
Conversely, if $({\overline{6}}_{\stt}) $ holds, we  know already that  ${(\stt)}_{\widetilde{\iota}} = w_{\widetilde{D}} = t_{\widetilde{D}}$ (by the fact that $({\overline{6}}_{\stt}) \Rightarrow ({{6}^\prime}_{\stt}) $)  and that  $ {(\stt)_{\overline{\iota}}} ={((\stt)_{\widetilde{\iota}})^j} =  (t_{\widetilde{D}})^j= ((v_{\widetilde{D}})^j)_{\!_f} $.

(b)   Let $D$ be a $\stf$-quasi Pr\"ufer domain.  If $\overline{D}$ 
is  $t$-linked to $(D, \star)$  then  
$\overline{D}$ is  a P${v_{\overline{D}}}$MD, 
since in this case $\overline{D} =\widetilde{D}$ 
(Lemma \ref{t-*-linked}  and Theorem \ref{um*}). 
 On the other hand, if $\overline{D}$ is not $t$-linked to $(D, \star)$,  
 then $\overline{D}$ is  a   P${{(\stt)}_{\overline{\iota}}}$MD (by (a))  and, since  in this case ${(\stt)}_{\overline{\iota}}$ is not a (semi)star operation on $\overline{D}$,  we may not deduce that $\overline{D}$ is a P${v_{\overline{D}}}$MD.   However, in the previous Example \ref{example}, even if $\overline{D}$ is not $t$-linked to $(D, t_D)$, we do have that   $\overline{D}$ is a  P${v_{\overline{D}}}$MD because for $\overline{H} \in \Max^{t_{\overline{D}}}(\overline{D})$ such that  $ \overline{H} \not\in \QMax^{(w_D)_{\overline{\iota}}}(\overline{D})$ we still have that $\overline{D}_{\overline{H}}$ is a valuation domain.

  (c)  Note that, if we replace $\stf$ with $\star$ in the conditions $(4_{\stf})$, $(8_{\stf})$, $(9_{\stf})$$(10_{\stf})$  and  $(14_{\stf})$  stated in Theorem \ref{um*}, we obtain:
     
     \begin{enumerate} \it 
     \item[$(4_{\star})$]  $ D \subseteq K$ is a   $\star$-primitive   extension.

\item[$(8_{\star})$] Every prime ideal of $\Na(D, \star)$   is extended from $D$.
 \item[$(9_{\star})$]   $\Na(D, \star)$    is a quasi-Pr\"ufer domain.
 \item[$(10_{\star})$] The integral closure of   $\Na(D, \star)$   is a Pr\"ufer domain. 
   \item[$(14_{\star})$]   $D$ is a   $t$-quasi-Pr\"ufer domain (or a UM\ \!$t$-domain)  and $\widetilde{\star} = w$.  

 \end{enumerate}
 It is trivial from the definitions that \it the previous conditions coincide with the analogous conditions stated for $\stf$ in Theorem \ref{um*}. \rm 
 
A natural question arises from this observation: is it possible to find suitable characterizations of the $\star$-quasi-Pr\"ufer domains, by ``weakening'' the remaining conditions in  Theorem \ref{um*}~\!?

 (d) Recall that Houston and Zafrullah have recently introduced the \it UM$v$-domains, \rm   i.e., the integral domains $D$, such that each upper to zero is a    maximal $v_{D[X]}$-ideal  of $D[X]$. It is known that UM$v$-domains are characterized by the fact that, for each upper to zero $P$, $\co_D(P)^{v_D} = D$ and $((P:P) = ) \ D[X] \subsetneq P^{-1}$ \cite[Theorem 2.2]{hz2}.
On the other hand, if $D$ is a UM$v$-domain and if $P$ is a $v_D$-prime ideal of $D$, then the integral closure of $D_P$ is a Pr\"ufer domain   \cite[Theorem 3.6]{hz2}, i.e., $D_P$ is a quasi-Pr\"ufer domain by Theorem \ref{qp}   ((1)$\Leftrightarrow$(6)). Therefore, by Lemma \ref{*qp}  ((iv)$\Rightarrow$(i)), a UM$v$-domain is a $v$-quasi-Pr\"ufer domain. 

Note also that a UM$v$-domain is not necessarily a $t$-quasi-Pr\"ufer domain (= UM$t$-domain). To see this,  let $D$ be a $v$-domain  (i.e.,  an integral domain such that each nonzero finitely generated ideal is $v$-invertible \cite[Theorem 34.6]{gilmer})     which is not a P$v$MD   (cf. \cite[Exercise 5, page 425]{gilmer} and also \cite[\S3]{d41}, \cite{h81} and \cite{ho72}).  A ring of this type must admit an upper to zero which is a maximal $v$-ideal but not a maximal $t$-ideal, since it is an integrally closed UM$v$ domain which is not a UM$t$-domain  (Remark \ref{umt} (d) and  \cite[Theorem 3.3 ((1)$\Leftrightarrow$(2))]{hz2}).  This  example  also shows that a $v$-quasi-Pr\"ufer domain need not be a 
$v_{_{\! f}}$-quasi-Pr\"ufer domain (cf. Example 2.2). 

\noindent \bf Question: \rm Is it possible to find a $v$-quasi-Pr\"ufer domain which is not a UM$v$-domain~\!?  

   (e)  Houston and Zafrullah \cite[Proposition 4.6]{hz2} proved that  $D$
is a UM$t$-domain if and only if each upper to zero of the form $(aX+b)K[X] \cap D[X]$, where $0 \neq a, b \in D$, is a maximal
$t$-ideal of $D[X]$ or, equivalently, each upper to zero of the form $(aX+b)K[X] \cap D[X]$, where $0 \neq a, b \in D$,  contains a   nonzero   polynomial $g$ with $\co_D(g)^t = D$  \cite[Theorem 1.4]{hz}.

A similar characterization holds for $\stf$-quasi-Pr\"ufer domains. More precisely, \it given a semistar operation $\star$ on an integral domain $D$, the following are equivalent:
\begin{enumerate}

\item[$(1_{\stf})$] $D$ is a $\stf$-quasi-Pr\"ufer domain.

\item[$(2^{\prime}_{\stf})$] 
Each upper to zero in $D[X]$ of the form $(aX+b)K[X] \cap D[X]$ 
contains a nonzero  polynomial $g$ with $\co_D(g)^\star = D^\star$. 

\item[$(2^{''}_{\stf})$] For each nonzero $h \in D[X]$, there exists $  0\neq g  \in hK[X]
\cap D[X]$ with $\co_D(g)^\star = D^\star$.
\end{enumerate} \rm 

$(1_{\stf})\!\Leftrightarrow\!(2^{\prime}_{\stf})$ By using the equivalence $(1_{\stf})\!\Leftrightarrow\!(4_{\stf})$ of Theorem \ref{um*} and the previous point (a), it is enough to show that
$Q := (aX+b)K[X] \cap D[X]$
 contains a  nonzero   polynomial $g$ with $\co_D(g)^\star = D^\star$ if and only if $u= -\frac{b}{a}$ is $\star$-primitive over $D$.

For the ``only if'' part, let $0 \neq g \in Q$ such that $\co_D(g)^\star = D^\star $.  Clearly  $g = (aX+b)h$, for some $h\in K[X]$. Then   $g(u) = \left(a\left(-\frac{b}{a}\right)+b\right)h(u) = 0$,     thus $u$ is $\star$-primitive over $D$. 
For the ``if''  part, suppose that $u \ (= -\frac{b}{a})$ is $\star$-primitive over $D$. Then there exists a nonzero polynomial 
$g \in D[X]$ such that $\co_D(g)^\star = D$ and $g(u) = 0$. Therefore in  $K[X]$ we have $g = (aX+b)h +r$, 
where $h \in K[X]$ and $r $ is a constant in $ K$.    Since $g(u)=0$, we have $r = 0$, and thus $g \in Q=(aX+b)K[X] \cap D[X]$.

The implication $(2^{''}_{\stf})\!\Rightarrow\!(2_{\stf})$ is obvious.

$(2_{\stf})\!\Rightarrow\! (2^{''}_{\stf})$  Let  $hK[X] ={\ell_1}{\ell_2}\dots {\ell_n}K[X]$, where $ \ell_i \in D[X]$  is irreducible in $K[X]$, for $1\leq i\leq n$.  Since $Q_i:= {\ell_i}K[X] \cap D[X]$ is an upper to zero, then we can find $0\neq g_i \in Q_i$ such that  $\co_D(g_i)^\star = D^\star$.   Then $g:=g_1g_2 \dots  g_n \in hK[X] \cap D[X]$ and it is not difficult to see that $\co_D(g)^\star = D^\star$. 

 (f) Note that,  from the equivalence $(1_{\stf})\Leftrightarrow(6_{\stt})$ in Theorem \ref{um*}   (or,  from   Corollary \ref{p*md-um*}), \ec we deduce that if  $\star$ is a (semi)star operation on $D$, then  $D$ is an integrally closed $\stf$-quasi-Pr\"ufer domain if and only if $D$ is a P$\stt$MD (or, equivalently, a P$\star$MD).  This result generalizes the statement on P$v$MDs recalled in Remark \ref{umt}  (d).
 
  \end{remark}

   \begin{corollary}   With the notation of Theorem \ref{um*}, we have that $(1_{\stf})$ is equivalent to
   
   \begin{enumerate}
\item [$(12^{\prime}_{\stf})$] Each $t$-linked overring $R$   to  $(D, \stf)$ is a   $t_R$-quasi-Pr\"ufer domain   and each $(\stt)_\iota$-maximal ideal   of $R$ is a $t_R$-ideal, where $\iota: D \hookrightarrow R$ is the canonical embedding.
\end{enumerate}
\end{corollary}

\begin{proof}  $(1_{\stf})\Rightarrow(12^{\prime}_{\stf})$  Note that from the proof $(1_{\stf})\Rightarrow(12_{\stf})$   of the previous Theorem 
\ref{um*} , we deduce, without assuming that $\star$ is a (semi)star operation on $D$,  that $R$ is a $(\widetilde{\star})_\iota$-quasi-Pr\"ufer domain. Henceforth $R$ is also  a  $t_R$-quasi-Pr\"ufer domain since $(\widetilde{\star})_\iota$ is a (semi)star operation of finite type on $R$. Now applying the implication $(1_{\stf})\Rightarrow(12_{\stf})$ to $R$ and to the (semi)star operation  $(\widetilde{\star})_\iota$, since $R$ is trivially $t$-linked to $(R, (\widetilde{\star})_\iota)$, we have in particular that each $(\stt)_\iota$-maximal ideal   of $R$ is a $t_R$-ideal.

$(12^{\prime}_{\stf})\Rightarrow (11_{\stf})$  If $P \in \QMax^{\stf}(D) = \QMax^{\stt}(D)$, then clearly
 $(D_P)^{\stt} = \bigcap \{D_PD_M \mid M \in \QMax^{\stf}(D)\} =D_P$ and so $D_P$ is $t$-linked to $(D, \stf)$  (Lemma \ref{t-*-linked} ((v)$\Rightarrow$(i))). Therefore, by assumption,  $D_P$ is a $t_{D_P}$-quasi-Pr\"ufer domain.  Moreover, clearly $PD_P$ is a maximal $(\stt)_\iota$-ideal of $D_P$ and so it is a $t_{D_P}$-ideal of $D_P$. Then $D_P$ is a quasi-Pr\"ufer domain by Corollary \ref{qp-umt} ((13)$\Rightarrow$(1)).
\end{proof}

\begin{corollary} \label{dim-umt}
If $D$ is a  $\stf$-quasi-Pr\"ufer domain,   then
\begin{enumerate}
\item[(a)]  If $P$ is a nonzero prime ideal of $D$ and if $P^{\stf} \neq D^\star$ (e.g., if $P$ is a   quasi-$\stf$-prime ideal of $D$), then   $P= P^{\stf} = P^t $. 
\item[(b)]   $\dim^{\stt}(D) = \dim^{\stf}(D) = \dim^t(D) =\dim^t(\Na(D, \star))= \dim(\Na(D, \star))$. 
\end{enumerate}
\end{corollary}

\begin{proof}
(a) It suffices to show that $P$ is a $t$-ideal. Let $Q$ be a
 quasi-$\stf$-maximal ideal of $D$ containing $P$. Then $D_Q$ is a
quasi-Pr\"ufer domain (Theorem \ref{um*} ($(1_{\stf})\Rightarrow(11_{\stf})$)), and since $PD_Q$ is a proper
prime ideal of $D_Q$, $PD_Q$ is a   prime $t$-ideal  in $D_Q$ (Corollary \ref{qp-umt}),  and hence $P = PD_Q \cap D$, is a
$t$-ideal of $D$  \cite[Lemma 3.17 (1)]{kang}.

(b) Note that  $\dim^{\stf}(D) = \dim^t(D)$  by (1) and $\dim^t(\Na(D, \star))= \dim(\Na(D, \star))$ by Corollary \ref{qp-umt} and 
Theorem \ref{um*}  ($(1_{\stf})\Rightarrow(9_{\stf})$).  
 Recall that $M \in \Max(\Na(D, \star))$ if and only if $M\cap D \in \QMax^{\stf}(D)$ \cite[Proposition 3.1 (5)]{fl}.  Since each prime ideal of  $\Na(D, \star)$  is extended from $D$ (Theorem \ref{um*}  ($(1_{\stf})\Rightarrow(8_{\stf})$)), we have  $\dim^{\stf}(D) =\dim(\Na(D, \star))$.   The first equality follows from the fact that the notions of $\stf$-quasi-Pr\"ufer domain and  $\stt$-quasi-Pr\"ufer domain 
coincide (Corollary \ref{<} (c))  and from the fact that $\Na(D, \stt) = \Na(D, \stf)= \Na(D, \star)$.   \end{proof}

It is well known that if $M$ is a maximal $t$-ideal of $D[X]$,
then either $M \cap D = (0)$ or $M = (M \cap D)[X]$
\cite[Proposition 1.1]{hz} and $I$ is a $t$-ideal of $D$
 if and only if  $I[X]$ is a $t$-ideal of $D[X]$ \cite[Corollary
2.3]{kang}. Thus  $\dim^t(D) \leq \dim^t(D[X]) \leq
2\dim^t(D))$  (cf. also \cite[page 169]{houston} and  \cite[Section 3]{wang}).

  Recall that Kang proved that if $D$ is a  P$v$MD then $\dim^t(D) = \dim(\Na(D, v)) $ \ec  \cite[Theorem 3.22]{kang}.
The following corollary extends Kang's result to the UM$t$-domains.

\begin{corollary}
Let $D$ be a  UM\ \!$t$-domain   which is not a field  and let $X$ be an indeterminate over $D$.
Then  $\dim^w(D) =  \dim^t(D) = \dim^t(D[X] )=\dim^t(\Na(D, v))= \dim(\Na(D, v))$.    
\end{corollary}

\begin{proof}
As we already remarked in general $\dim^t(D) \leq  \dim^t(D[X])$. 
Let $Q$ be a maximal $t$-ideal of $D[{X}]$. If $Q \cap D = (0)$, then obviously  $\hgt(Q)= 1 \leq \dim(\Na(D,v))$.  If $Q \cap
D \neq (0)$, then $Q = (Q\cap D)[X]$ and hence $Q \cap \mathcal{N}^v = \emptyset$. Therefore  $Q\Na(D,v) \neq \Na(D,v)$ and so  $\hgt(Q) \leq
\dim(\Na(D, v))$, hence 
 $ \dim^t(D[X] )\leq \dim(\Na(D, v))$.  
The conclusion follows easily from Corollary \ref{dim-umt} (b).
\end{proof}

\begin{remark} \rm  (a) Note that, for a UM$t$-domain,    Wang  \cite[Theorem 2.6]{wang2} proved already the equality $\dim^w(D) = \dim(\Na(D, v))$.  

(b) It is clear  that, in general, $\dim^t(D) \leq \dim^w(D)$, since each $t$-ideal is also a $w$-ideal and it is easy to see that (in the non UM$t$-domain case) it can happen that $\dim^t(D) \neq  \dim^w(D)$. For instance, let $R$ be a quasi-local factorial domain of dimension $ 3$ with maximal ideal $M$. Set 
$F:= R/M$ and let $\varphi : R \rightarrow F$ be the
canonical homomorphism. Assume that   $k$ is a proper subfield of $F$, set $D := \varphi^{-1}(k)$ and let  $Q$ is a prime ideal of $D$ and $R$  such that
 $\hgt(Q) = 2$.  Clearly, since $R$ is a UFD and $M = (D:R)$, then   $ M = M^{v_D} = M^{w_D} \subsetneq M^{w_R} = M^{t_R} = M^{v_R} = R$ and 
$Q = Q^{w_D} \subsetneq Q^{t_R} = Q^{v_R} = R$  (note that $Q = Q^{w_D}$\!, since $\Max^{w_D}(D) = \Max(D)$ and so $w_D$ coincides with the identity (semi)star operation on $D$).   Let $I \subseteq Q$ be a nonzero finitely generated ideal of   $D$ with $(R:(R:I)) = (R:(R:IR))  = R$ or, equivalently, $(R:I) =R$.  
Hence $(D:I) \subseteq (R:I) = R = (M:M) =(D:M) \subseteq (D:I)$ and so  $(D:I) = R$. Therefore $I^{v_D} =(D:(D:I)) = (D:R) =M$ and so  $M = I^{v_D} \subseteq Q^{t_D} \subseteq M^{t_D} = M$. Henceforth  
$Q^{t_D} = M$.   Therefore we have
$\dim^t(D) = 2 \lneq \dim(D) = \dim^w(D)=3$.

(c)    It is well known that an integral domain $D$ is  Pr\"ufer domain (resp., P$v$MD),
if and only if each nonzero two generated ideal of $D$ is
invertible (resp., $t$-invertible) \cite[Theorem
22.1]{gilmer}  (resp., \cite[Lemma 1.7]{mmz}).     In case $\ast$ is a star operation  of finite type,     it is known that $D$ is  P$\ast$MD if and only if each (nonzero) two generated ideal of $D$ is $\ast$-invertible  \cite[Theorem 1.1]{hmm}.   It is natural to ask whether a similar result holds in the semistar setting.
 Let $\star$ be a semistar operation on an integral domain $D$. Recall that, in \cite[Theorem 2.3]{fo-pi},  it is shown that for $I \in \boldsymbol{f}(D)$,  $I$ is $\stf$-invertible if and only if $ ID_Q$ is principal, for each $P \in \QMax^{\stf}(D)$.  Moreover, it is well known that, for a local domain, the following properties are equivalent  \cite[Theorem 22.1]{gilmer}:
 \begin{enumerate}
 \item[(i)]  Every nonzero finitely generated ideal is principal;
 \item[(ii)]  Every two generated is principal;
 \item[(iii)]  $R$ is a valuation domain.
 \end{enumerate}
 On the other hand, $D$ is a P$\star$MD if and only if $D_P$ is a valuation domain,  for each $P \in \QMax^{\stf}(D)$ \cite[Theorem 3.1]{fjs}. Therefore, by the previous considerations it follows that \it $D$ is a P$\star$MD if and only if each (nonzero) two generated ideal of $D$ is    $\star_f$-invertible.   \rm 


 \end{remark}

\begin{center}
\sc
 Acknowledgments \rm
\end{center}
 \rm  During the preparation of this paper, the second named author
was partially supported by a grant PRIN-MiUR.


\end{document}